\documentclass[a4paper,10pt,DIV=12]{scrartcl}

\usepackage{amsmath,amssymb,amsthm}
\usepackage[T1]{fontenc}
\usepackage{lmodern}
\usepackage{bbm}            
\usepackage{graphicx}
\usepackage{float}
\usepackage[dvipsnames]{xcolor} 
\usepackage{tikz,pgf}
\usepackage{booktabs}
\usepackage[colorlinks=true,linkcolor=RoyalBlue,citecolor=PineGreen,urlcolor=RoyalBlue]{hyperref}
\usepackage{caption}
\usepackage{subcaption}

\newcommand{\lam}{\mathrm{Lam}}
\newcommand{\lamII}{\lam_2}
\newcommand{\lamIII}{\lam_3}
\newcommand{\Mtwo}{M_2}
\newcommand{\Mthree}{M_3}
\newcommand{\maxlamII}[1]{\Mtwo(#1)}
\newcommand{\maxlamIII}[1]{\Mthree(#1)}
\newcommand{\maxlamIIn}{\maxlamII{n}}
\newcommand{\maxlamIIIn}{\maxlamIII{n}}
\newcommand{\modop}{\,\mathrm{mod}\,}

\renewcommand{\leq}{\leqslant}
\renewcommand{\geq}{\geqslant}
\newcommand{\N}{\mathbbm N}
\newcommand{\Z}{\mathbbm Z}

\newcommand{\R}{\mathbbm R}
\newcommand{\C}{\mathbbm C}

\newtheorem{theorem}{Theorem}

\newtheorem{definition}[theorem]{Definition}
\theoremstyle{definition}

\newtheorem{oproblem}{Open Problem}
\newtheorem{conjecture}[oproblem]{Conjecture}

\colorlet{ncol}{green!40!black}
\tikzstyle{vertex}=[circle, draw, fill=black, inner sep=0pt, minimum size=4pt]
\tikzstyle{nvertex}=[circle, draw=ncol, fill=ncol, inner sep=0pt, minimum size=4pt]
\tikzstyle{edge}=[thick]
\tikzstyle{edgeq}=[gray!60,densely dashed,thick]
\tikzstyle{nedge}=[ncol,thick]
\tikzstyle{oedge}=[red!60!black,thick]
\newcommand{\vertex}{\node[vertex]}

\begin{document}

\title{Lower bounds on the number of realizations of rigid graphs}

\author{%
  Georg Grasegger%
  \thanks{Johann Radon Institute for Computational and Applied
    Mathematics (RICAM), Austrian Academy of Sciences,
    Altenberger Stra\ss e 69, 4040 Linz, Austria} \and Christoph
  Koutschan\footnotemark[1]%
  \and Elias Tsigaridas%
  \thanks{ Sorbonne Universit{\'e}s, \textsc{UPMC} Univ Paris 06,
    \textsc{CNRS}, \textsc{INRIA}, {Laboratoire d'Informatique de
      Paris 6 (\textsc{LIP6}), {\'E}quipe \textsc{PolSys}, 4 place
      Jussieu, 75252 Paris Cedex 05, France}}
}

\date{March 27, 2018\footnote{This paper has been published in Experimental
    Mathematics, 2018, DOI:
    \href{http://dx.doi.org/10.1080/10586458.2018.1437851}{10.1080/10586458.2018.1437851}.}}

\maketitle

\begin{abstract}
  Computing the number of realizations of a minimally rigid graph is a
  notoriously difficult problem.
  Towards this goal, for graphs that are minimally rigid in the
  plane, we take advantage of a recently published algorithm, which is the
  fastest available method, although its complexity is still exponential.
  Combining computational results with the theory of constructing new
  rigid graphs by gluing, we give a new lower bound on the maximal
  possible number of (complex) realizations for graphs with a given number of
  vertices.  We extend these ideas to rigid graphs in three
  dimensions and we derive similar lower bounds, by exploiting data
  from extensive Gr\"obner basis computations.
\end{abstract}

\section{Introduction}

The theory of rigid graphs forms a fascinating research area in the
intersection of graph theory, computational (algebraic) geometry, and
algorithms. Besides being a very interesting mathematical subject,
rigid graphs and the underlying theory of Euclidean
distance geometry
have a huge number of applications ranging from robotics
\cite{FL94,WalHus-9bar-07,WalHus-9bar-3d-07} and bioinformatics
\cite{EmiMou99,JRKT01,LiLaMuMa-molecular-11,LiMaLeLaMu-rigid-14,MuLaLiMa-book-12} to sensor network
localization \cite{ZhSoYe-urigid-10} and architecture
\cite{Emmerich-rigid-88}.  Upper and lower bounds on
the number of realizations (embeddings) of rigid graphs are of
great importance as they quantify the difficulty of the problem(s) at
hand that we are interested in.

We first give some definitions, to set the context of our study.  Let
$G$ be a graph and provide to~$\R^d$ the Euclidean metric; in this way we
obtain the Euclidean $d$-dimensional  space.  By specifying the
coordinates of the vertices of~$G$ in~$\R^d$ we obtain a \emph{realization},
or \emph{embedding}, of~$G$ in~$\R^d$.
If there is no
continuous deformation of the graph that preserves the edge lengths, then the embedding is called \emph{rigid}.
A graph~$G$ is said to be \emph{generically rigid} in~$\R^d$ if and only if all of its generic realizations are rigid. 
In the case of~$\R^2$ these graphs are also known as \emph{Laman graphs}.

Given a generically rigid graph in~$\R^d$, together with generic edge lengths,
we can embed it in the Euclidean $d$-space in a finite number of
ways, modulo rigid motions (translations and rotations).
It is of great interest to provide tight bounds for the number of
embeddings of such graphs, modulo rigid motions, for any~$d$.
Our results provide lower bounds for $d =2$ and $d = 3$.

\subsection{Previous work}
\label{sec:previous-work}

The first bounds on the number of realizations of rigid graphs, using
degree bounds from algebraic geometry, are due to Borcea and Streinu
\cite{Borcea2004}.  They rely on the theory of distance matrices and
on bounds of determinantal varieties. This results in the upper bounds
$\binom{2n-4}{n-2} =\Theta( 4^{n}/ \sqrt{n})$ for graphs in~2D, and
$\frac{2^{n-3}}{n-2}\binom{2n-6}{n-3} = \Theta( 8^{n}/(n\sqrt{n}))$
for graphs in~3D, where $n$ denotes the number of vertices.
Steffens and Theobald \cite{SteTheo-cg-10} improved
these bounds by exploiting the sparsity of the underlying polynomial
systems.  These bounds were further improved by applying additional
tricks to take advantage of the sparsity and the common
sub-expressions that appear in the polynomial systems~\cite{Emiris2013,Emiris2009}.
A direct application of the mixed
volume techniques, which roughly speaking capture the sparsity of a
polynomial system, yield a bound of $4^{n-2}$ for the planar case.
If we also take into account the degree of the vertices, then in the 2D case, for a Laman graph with $k \geq 4$
degree-2 vertices, the number of planar embeddings of $G$ is bounded from
above by $2^{k-4} 4^{n-k}$.  For the 3D case, when the graph is the
1-skeleton of a simplicial polyhedron with $k \geq 9$ degree-3
vertices, then the number of embeddings is bounded from above by
$2^{k-9} 8^{n-k}$.

The state-of-the-art result is the recent paper~\cite{SCgroup}
that provides an algorithm for computing the number of complex realizations of Laman graphs.
The algorithm recursively computes these numbers by lifting the problem to pairs of graphs.
Arguments from tropical geometry are used to show the correctness of the algorithm,
while the computations themselves are then purely combinatorial (an implementation can be found at \cite{SCgroupElectronicMaterial}).
With help of the algorithm, the number of realizations of all Laman graphs up to 12 vertices were computed.
We exploit these data in the present paper.

The first lower bounds for graphs in 2D were 
$24^{\lfloor (n-2)/4 \rfloor}$ (approx.\ $2.21^n$) and
$2\cdot 12^{\lfloor (n-3)/3 \rfloor}$ (approx.\ $2.29^n$), that exploited a
gluing process using a
caterpillar, resp.\ fan construction \cite{Borcea2004}, see also
\cite{Emiris2009}.  Both constructions use the three-prism graph (sometimes also called Desargues graph) as a
building block, which is a graph with 6 vertices and 24 embeddings.  More recent
lower bounds are $2.30^n$ from \cite{EmirisMoroz} and $2.41^n$
from~\cite{Jackson2018}.
For graphs in 3D, the only known lower bound is
$ 16^{\lfloor (n-3)/3\rfloor}$ (approx.\ $2.52^n$) for $n\geq 9$, which uses
a cyclohexane caterpillar as building block \cite{Emiris2009}.

\subsection{Our contribution}
\label{sec:our-contrib}

We present lower bounds on the maximal number of planar, resp.\ spatial,
embeddings (up to rigid motions) of minimally rigid graphs with a prescribed
number of vertices. However, we relax the condition that the
embeddings take place in~$\R^d$.  Instead, we compute the number of
\emph{complex} Euclidean embeddings, that is embeddings in~$\C^d$.  
In this complex setting, even the edge lengths may be assumed to be
complex numbers. Clearly, the number of complex embeddings is an upper
bound on the number of real embeddings.

Using the novel algorithm developed in~\cite{SCgroup} we compute
the exact number of planar embeddings
for graphs with a relatively small number of vertices.
In contrast, the number of spatial embeddings is computed probabilistically
by means of Gr\"obner bases. Then we introduce
techniques to ``glue'' an arbitrary number of such small graphs in order to
produce graphs with a high number of vertices (and edges) that
preserve rigidity.  The gluing process (see
Sections~\ref{sec:constr-2d} and \ref{sec:constr-3d}) allows us to derive
the number of embeddings of the final graph from the number of
embeddings of its components, and in this way we derive a lower bound
for the number of embeddings in~$\C^2$ (Theorem~\ref{thm:lower-bd-2d})
and in~$\C^3$ (Theorem~\ref{thm:lower-bd-3d}).  We emphasize that the
gluing techniques are quite general and can be extended to arbitrary
dimensions. Moreover, to identify those small graphs that realize
the maximum number of embeddings and that can be the building blocks for the
gluing process, we perform extensive experiments. We use the
state-of-the-art computer algebra tools to count the number of
embeddings as the maximum number of complex solutions of polynomial
systems.

If we were able to compute the number of embeddings of the small
graphs in~$\R^d$, for example by using the approach proposed
in~\cite{EmirisMoroz}, then we could transfer our lower bounds on the
complex embeddings to the number of real embeddings, by applying the
very same gluing process; see also \cite{Jackson2018} for gluing
processes using the caterpillar graph. There it is also hinted
that the numbers of real and
complex embeddings do not match in general. It is a very interesting
problem to quantify this gap. On the one hand, one can construct
infinite families of graphs for which the ratio between real and
complex embeddings tends to zero. On the other hand, there are graphs,
see~\cite{EmirisMoroz} for a nontrivial example, where edge lengths
can be found such that there exist as many real embeddings as complex
ones.

\subsection{Organization of the paper}
The paper is structured as follows: First (Section~\ref{sec:dim-2}) we present
the construction of the lower bounds for the planar case, and in
Section~\ref{sec:dim-3} we present the lower bounds for the spatial case. In
Section~\ref{sec:constr-2d} we describe three constructions (gluing processes)
for producing infinite families of rigid graphs. Then, in
Section~\ref{sec:graphs-2d}, we discuss several strategies to identify
expedient graphs that are suitable for these constructions. They lead to new
lower bounds, which is discussed in Section~\ref{sec:bounds-2d}.

Throughout the paper we represent a graph by the integer obtained by
flattening the upper triangular part of its adjacency matrix and interpreting
this binary sequence as an integer.  For further details we refer to
Appendix~\ref{appendix:encoding}. There we also collect the encodings
of all graphs mentioned throughout the paper.

\section{Dimension 2}
\label{sec:dim-2}

We begin our study with the case of planar embeddings. For this purpose,
we recall the definitions of some fundamental notions. The goal in this
section is to derive lower bounds for the quantity $\maxlamIIn$, introduced
in Definition~\ref{def:laman-number-2d} below.

\begin{definition}\label{definition:laman_graph}
	A \emph{Laman graph}~\cite{Laman1970} is a graph $G=(V,E)$ such that $\vert E\vert = 2\vert
	V\vert-3$, and such that $\vert E' \vert \leq 2\vert V' \vert-3$ holds for
	every subgraph $G'=(V',E')$ of~$G$.
\end{definition}

\begin{definition}\label{def:laman-number-2d}
	For a Laman graph $G=(V,E)$ we define $\lamII(G)$, called the \emph{Laman number} of~$G$,
	to be the number of (complex) planar embeddings that a generic labeling $\lambda\colon E\to\C$
        (the ``edge lengths'' of~$G$) admits. Moreover, we define $\maxlamIIn$ to be the largest Laman
	number that is achieved among all Laman graphs with $n$ vertices.
\end{definition}

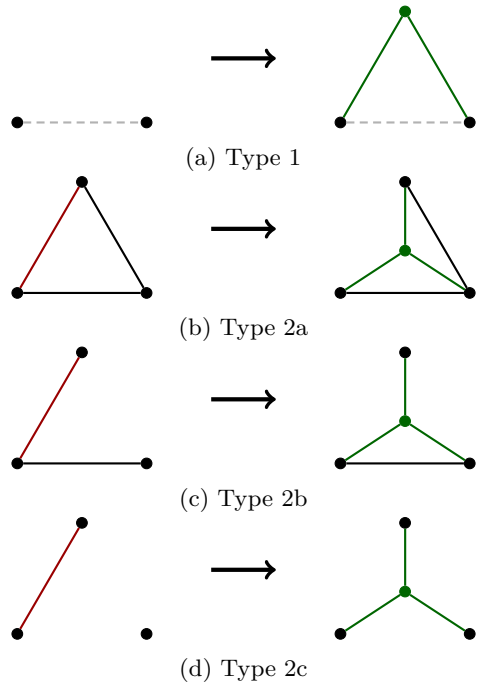
\begin{figure}[H]
  \begin{center}
		\begin{subfigure}[b]{0.98\textwidth}
      \begin{center}
				\begin{tikzpicture}[scale=1.7]
					\vertex (a) at (0,0) {};
					\vertex (b) at (1,0) {};
					\vertex (d) at (2.5,0) {};
					\vertex (e) at (3.5,0) {};
					\node[nvertex] (f) at (3,0.866025) {};
					\draw[edgeq] (a)edge(b) (d)edge(e);
					\draw[nedge] (d)edge(f) (e)edge(f);
					\draw[ultra thick,->] (1.5,0.5) -- (2,0.5);
				\end{tikzpicture}
				\caption{Type~1}\label{fig:type1}
      \end{center}
    \end{subfigure}
		\\
		\begin{subfigure}[b]{0.98\textwidth}
      \begin{center}
				\begin{tikzpicture}[scale=1.7]
					\vertex (a) at (0,0) {};
					\vertex (b) at (1,0) {};
					\vertex (c) at (0.5,0.866025) {};
					\vertex (d) at (2.5,0) {};
					\vertex (e) at (3.5,0) {};
					\vertex (f) at (3,0.866025) {};
					\node[nvertex] (g) at (3,0.33) {};
					\draw[edge] (a)edge(b) (b)edge(c) (d)edge(e) (e)edge(f);
					\draw[oedge] (a)edge(c);
					\draw[nedge]  (d)edge(g) (e)edge(g) (f)edge(g);
					\draw[ultra thick,->] (1.5,0.5) -- (2,0.5);
				\end{tikzpicture}
				\caption{Type~2a}\label{fig:type2a}
      \end{center}
    \end{subfigure}
    \\
    \begin{subfigure}[b]{0.98\textwidth}
      \begin{center}
				\begin{tikzpicture}[scale=1.7]
					\vertex (a) at (0,0) {};
					\vertex (b) at (1,0) {};
					\vertex (c) at (0.5,0.866025) {};
					\vertex (d) at (2.5,0) {};
					\vertex (e) at (3.5,0) {};
					\vertex (f) at (3,0.866025) {};
					\node[nvertex] (g) at (3,0.33) {};
					\draw[edge] (a)edge(b)  (d)edge(e);
					\draw[oedge] (a)edge(c);
					\draw[nedge]  (d)edge(g) (e)edge(g) (f)edge(g);
					\draw[ultra thick,->] (1.5,0.5) -- (2,0.5);
				\end{tikzpicture}
				\caption{Type~2b}\label{fig:type2b}
      \end{center}
    \end{subfigure}
    \\
    \begin{subfigure}[b]{0.98\textwidth}
      \begin{center}
				\begin{tikzpicture}[scale=1.7]
					\vertex (a) at (0,0) {};
					\vertex (b) at (1,0) {};
					\vertex (c) at (0.5,0.866025) {};
					\vertex (d) at (2.5,0) {};
					\vertex (e) at (3.5,0) {};
					\vertex (f) at (3,0.866025) {};
					\node[nvertex] (g) at (3,0.33) {};
					\draw[oedge] (a)edge(c);
					\draw[nedge]  (d)edge(g) (e)edge(g) (f)edge(g);
					\draw[ultra thick,->] (1.5,0.5) -- (2,0.5);
				\end{tikzpicture}
				\caption{Type~2c}\label{fig:type2c}
      \end{center}
    \end{subfigure}
  \end{center}
  \caption{Henneberg steps of different types in dimension~2; A dashed line indicates that this edge can exist but does not need to.}
  \label{fig:Henneberg}
\end{figure}

In \cite{Henneberg} Laman graphs are characterized to be constructible from a single edge
by a sequence of two types of steps (see Figure~\ref{fig:Henneberg}).
We call them Henneberg steps of type~1 and type~2 respectively.
The steps of type~2 can be further classified according to additional occurring edges.

\begin{figure}[H]
  \begin{center}
    \begin{tikzpicture}[scale=1.7]
      \begin{scope}[xshift=-0.5cm]
				\vertex (a1) at (0,0) {};
				\vertex (a2) at (-0.4,0) {};
				\vertex (b1) at (1,0) {};
				\vertex (b2) at (1.4,0) {};
				\vertex (c) at (0.5,0.866025) {};
				\vertex (d) at (0.5,0) {};
				\draw[dashed,black!70!white] (-0.2,0) circle[x radius=0.35cm, y radius=0.15cm];
				\draw[dashed,black!70!white] (1.2,0) circle[x radius=0.35cm, y radius=0.15cm];
				\node[black!70!white,align=center] at (-0.2,0) {\tiny\ldots};
				\node[black!70!white,align=center] at (1.2,0) {\tiny\ldots};
				\draw[edge] (a1)edge(c)  (a2)edge(c) (d)edge(c);
				\draw[oedge] (b1)edge(c) (b2)edge(c);
			\end{scope}
			\draw[ultra thick,->] (1.5,0.5) -- (2,0.5);
			\begin{scope}[xshift=3cm]
			  \vertex (a1) at (0,0) {};
				\vertex (a2) at (-0.4,0) {};
				\vertex (b1) at (1,0) {};
				\vertex (b2) at (1.4,0) {};
				\vertex (c) at (0.3,0.866025) {};
				\node[nvertex] (c2) at (0.7,0.866025) {};
				\vertex (d) at (0.5,0) {};
				\draw[dashed,black!70!white] (-0.2,0) circle[x radius=0.35cm, y radius=0.15cm];
				\draw[dashed,black!70!white] (1.2,0) circle[x radius=0.35cm, y radius=0.15cm];
				\node[black!70!white,align=center] at (-0.2,0) {\tiny\ldots};
				\node[black!70!white,align=center] at (1.2,0) {\tiny\ldots};
				\draw[edge] (a1)edge(c)  (a2)edge(c) (d)edge(c);
				\draw[nedge] (b1)edge(c2) (b2)edge(c2) (d)edge(c2) (c)edge(c2);
			\end{scope}
    \end{tikzpicture}
  \end{center}
  \caption{Vertex splitting}
  \label{fig:vertexsplitting}
\end{figure}
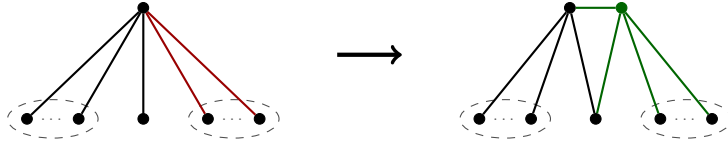

It is well known that a Henneberg step of type~1 always increases the Laman number by a factor of~$2$.
So far it is not known by which factor a Henneberg step of type~2 might increase the Laman number.
As mentioned in \cite{Jackson2018} there are Henneberg steps of type~2 which do increase the Laman number by a factor of less than~$2$.
Vertex splitting is another construction preserving rigidity (see Figure~\ref{fig:vertexsplitting}).
In \cite{Jackson2018} it is shown that vertex splitting increases the Laman number by a factor of at least two.
Henneberg steps of type 2a and 2b are special cases of vertex splitting.
Hence, only type 2c can yield a factor of less than two.
Table~\ref{table:Laman_Increase_2d} shows some increases of Laman numbers, given a certain Laman graph $G$ and constructing a new one $G'$ by a single Henneberg step.
\begin{table}[H]
  \begin{center}
    \begin{tabular}{lrrrrr}
      \toprule
      Type & $G$              & $\lamII(G)$ & $G'$              & $\lamII(G')$ & Factor\\
      \midrule
      2c   & 1269995          & 56  & 31004235             & 96   & 1.71  \\
      2c   & 7916             & 24  & 481867               & 44   & 1.83  \\
      2b   & 186013           & 32  & 170989214            & 136  & 4.25  \\
      2c   & 183548           & 32  & 170989214            & 136  & 4.25  \\
      2c   & 20042142         & 64  & 11177989553          & 344  & 5.37  \\
      2c   & 4593214614       & 128 & 22301628505804       & 808  & 6.31  \\
      2c   & 1248809223262    & 256 & 2960334732174949     & 1976 & 7.72  \\
      2c   & 1710909647295913 & 512 & 15006592507478215906 & 4816 & 9.41 \\
      \bottomrule
    \end{tabular}
  \end{center}
  \caption{Henneberg constructions and increase of Laman numbers}
  \label{table:Laman_Increase_2d}
\end{table}

\subsection{Constructions}
\label{sec:constr-2d}

We discuss different constructions of infinite families of Laman
graphs $(G_n)_{n\in\N}$ with $G_n$ having $n$ vertices. We do this in a way
such that we know precisely the Laman number for each member of the family.
This directly leads to a lower bound on $\maxlamIIn$.
The ideas of these constructions are described in~\cite{Borcea2004};
they were used to get lower bounds by connecting several three-prism graphs at a common basis.
Here, we generalize them in order to connect any Laman graphs at an arbitrary Laman base.
We present three such constructions.

\subsubsection{Caterpillar construction}
\label{sec:caterpillar}

The ``caterpillar construction''~\cite{Borcea2004} works as
follows: place $k$ copies of a Laman graph $G=(V,E)$ in a row and
connect every two neighboring ones by means of a shared edge (see
Figure~\ref{figure:caterpillar}).  Alternatively, one can let all $k$
graphs share the same edge. In any case, the resulting assembly has
$2+k(|V|-2)$ vertices and its Laman number is $\lamII(G)^k$, since each of
the $k$ copies of $G$ can achieve all its $\lamII(G)$ different embeddings,
independently of what happens with the other copies. Hence, among all
Laman graphs with $n=2+k(|V|-2)$ vertices there exists one with
$\lamII(G)^k$ embeddings. If the number of vertices $n$ is not of the form
$2+k(|V|-2)$ then we can use the previous caterpillar graph with
$\lfloor(n-2)/(|V|-2)\rfloor$ copies of $G$ and perform some
Henneberg steps of type~1; as we mentioned earlier, each of these
steps doubles the Laman number.  Summarizing, for any Laman graph~$G$,
we obtain the following lower bound from the caterpillar construction:
\begin{equation}\label{eq:bound_cater}
  \maxlamIIn \geq 2^{(n-2)\modop(|V|-2)} \cdot \lamII(G)^{\lfloor(n-2)/(|V|-2)\rfloor}\qquad (n\geq2).
\end{equation}
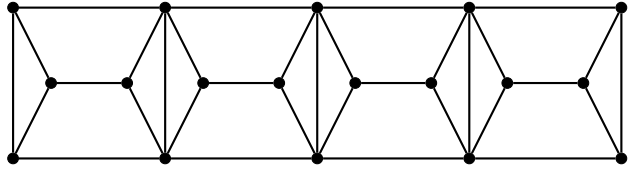
\begin{figure}[H]
	\begin{center}
		\begin{tikzpicture}[scale=1]
			\vertex (a) at (0.00,0.00) {};
			\vertex (b) at (0.00,2.00) {};
			\vertex (c) at (0.50,1.00) {};
			\vertex (d) at (1.50,1.00) {};
			\vertex (e) at (2.00,0.00) {};
			\vertex (f) at (2.00,2.00) {};
			\vertex (g) at (2.50,1.00) {};
			\vertex (h) at (3.50,1.00) {};
			\vertex (i) at (4.00,0.00) {};
			\vertex (j) at (4.00,2.00) {};
			\vertex (k) at (4.50,1.00) {};
			\vertex (l) at (5.50,1.00) {};
			\vertex (m) at (6.00,0.00) {};
			\vertex (n) at (6.00,2.00) {};
			\vertex (o) at (6.50,1.00) {};
			\vertex (p) at (7.50,1.00) {};
			\vertex (q) at (8.00,0.00) {};
			\vertex (r) at (8.00,2.00) {};
			\draw[edge] (a)edge(b) (a)edge(c) (b)edge(c) (a)edge(e) (b)edge(f)
				(c)edge(d) (d)edge(e) (e)edge(f) (d)edge(f) (e)edge(g) (f)edge(g)
				(e)edge(i) (f)edge(j) (g)edge(h) (h)edge(i) (i)edge(j) (h)edge(j)
				(i)edge(k) (j)edge(k) (i)edge(m) (j)edge(n) (k)edge(l) (l)edge(m)
				(m)edge(n) (l)edge(n) (m)edge(o) (n)edge(o) (m)edge(q) (n)edge(r)
				(o)edge(p) (p)edge(q) (q)edge(r) (p)edge(r);
		\end{tikzpicture}
	\end{center}
	\caption{Caterpillar construction with $4$ copies of the three-prism graph.}
	\label{figure:caterpillar}
\end{figure}

\subsubsection{Fan construction}
\label{sec:fan}

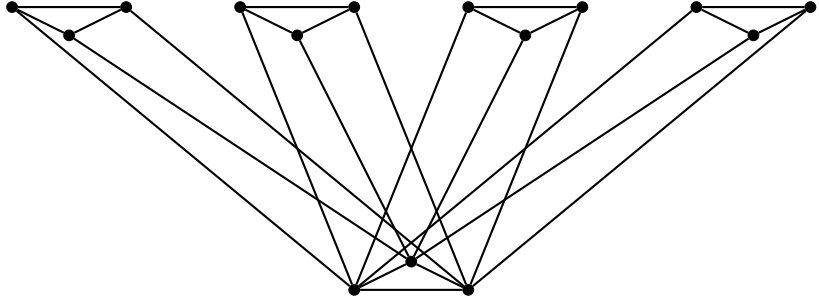
\begin{figure}[H]
	\begin{center}
		\begin{tikzpicture}[scale=0.75]
			\vertex (a) at (-1.00,0.00) {};
			\vertex (b) at (1.00,0.00) {};
			\vertex (c) at (0.00,0.50) {};
			\vertex (d) at (-7.00,5.00) {};
			\vertex (e) at (-5.00,5.00) {};
			\vertex (f) at (-6.00,4.50) {};
			\vertex (g) at (-3.00,5.00) {};
			\vertex (h) at (-1.00,5.00) {};
			\vertex (i) at (-2.00,4.50) {};
			\vertex (j) at (1.00,5.00) {};
			\vertex (k) at (3.00,5.00) {};
			\vertex (l) at (2.00,4.50) {};
			\vertex (m) at (5.00,5.00) {};
			\vertex (n) at (7.00,5.00) {};
			\vertex (o) at (6.00,4.50) {};
			\draw[edge] (a)edge(b) (b)edge(c) (c)edge(a) (d)edge(e) (e)edge(f)
				(f)edge(d) (d)edge(a) (e)edge(b) (f)edge(c) (g)edge(h) (h)edge(i)
				(i)edge(g) (g)edge(a) (h)edge(b) (i)edge(c) (j)edge(k) (k)edge(l)
				(l)edge(j) (j)edge(a) (k)edge(b) (l)edge(c) (m)edge(n) (n)edge(o)
				(o)edge(m) (m)edge(a) (n)edge(b) (o)edge(c);
		\end{tikzpicture}
	\end{center}
	\caption{Fan construction with $4$ copies of the three-prism graph.}
	\label{figure:fan}
\end{figure}

The second construction we employ is called ``fan construction'': take a
Laman graph $G=(V,E)$ that contains a triangle (i.e., a $3$-cycle), and glue $k$ copies of $G$
along that triangle (see Figure~\ref{figure:fan}). Once we fix one of the two
possible embeddings of that triangle, each copy of $G$ admits $\lamII(G)/2$
embeddings. The remaining $\lamII(G)/2$ embeddings are obtained by mirroring, i.e.,
by using the second embedding of the common triangle. Similarly as before, the
assembled fan is a Laman graph with $3+k(|V|-3)$ vertices that admits
$2\cdot(\lamII(G)/2)^k$ embeddings. Hence, we get the following lower bound:
\begin{equation}\label{eq:bound_fan}
  \maxlamIIn \geq 2^{(n-3)\modop(|V|-3)}\cdot 2\cdot\left(\frac{\lamII(G)}{2}\right)^{\!\lfloor(n-3)/(|V|-3)\rfloor} \qquad (n\geq3).
\end{equation}
While the caterpillar construction can be done with any Laman graph, this is
not the case with the fan. For example, the Laman graph with 12 vertices displayed in 
Figure~\ref{figure:max_graphs} has no $3$-cycle and therefore cannot be used for
the fan construction (see also Table~\ref{table:bounds}).

\subsubsection{Generalized fan construction}
\label{sec:gfan}

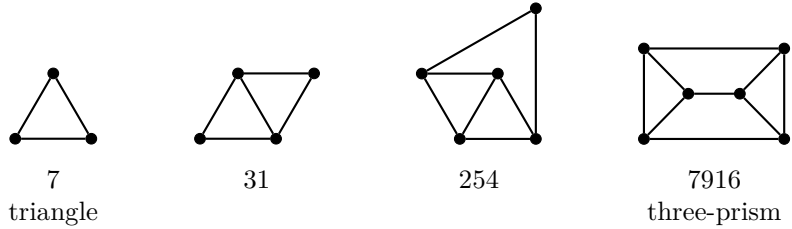
\begin{figure}[H]
\begin{center}
  \setlength{\tabcolsep}{18pt}
  \begin{tabular}{cccc}
    \begin{tikzpicture}
      \vertex (a) at (0,0) {};
			\vertex (b) at (1,0) {};
			\vertex (c) at (0.5,0.866025) {};
			\draw[edge] (a)edge(b) (b)edge(c) (a)edge(c);
    \end{tikzpicture}
    &
    \begin{tikzpicture}
      \vertex (a) at (0,0) {};
			\vertex (b) at (1,0) {};
			\vertex (c) at (0.5,0.866025) {};
			\vertex (d) at (1.5,0.866025) {};
			\draw[edge] (a)edge(b) (b)edge(c) (a)edge(c) (b)edge(d) (c)edge(d);
    \end{tikzpicture}
    &
    \begin{tikzpicture}[rotate=-60]
      \vertex (a) at (0,0) {};
			\vertex (b) at (1,0) {};
			\vertex (c) at (0.5,0.866025) {};
			\vertex (d) at (1.5,0.866025) {};
			\vertex (e) at (0,1.732050) {};
			\draw[edge] (a)edge(b) (b)edge(c) (a)edge(c) (b)edge(d) (c)edge(d) (a)edge(e) (d)edge(e);
    \end{tikzpicture}
    &
    \begin{tikzpicture}[yscale=0.04,xscale=0.05,rotate=90]
      \vertex (a) at (-15.,-18.4482) {};
			\vertex (b) at (15,-18.4482) {};
			\vertex (c) at (0.,-6.78132) {};
			\vertex (d) at (-15.,18.4482) {};
			\vertex (e) at (15.,18.4482) {};
			\vertex (f) at (0.,6.78132) {};
			\draw[edge]  (a)edge(b) (a)edge(c) (a)edge(d) (b)edge(c) (b)edge(e) (c)edge(f) (d)edge(e) (d)edge(f) (e)edge(f);
    \end{tikzpicture}
    \\[1ex]
    7 & 31 & 254 & 7916\\
    triangle & & & three-prism
  \end{tabular}
\end{center}
\caption{Bases for the generalized fan construction and their encodings.}
\label{figure:fanbases}
\end{figure}

As a third construction, we propose a generalization of the fan construction:
instead of a triangle, we may use any Laman subgraph $H=(W,F)$ of $G$ for
gluing. Using $k$ copies of $G$, we end up with a fan consisting of
$|W|+k(|V|-|W|)$ vertices and Laman number at least $\lamII(H)\cdot(\lamII(G)/\lamII(H))^k$. Here we
assume that the embeddings of $G$ are divided into $L(H)$ equivalence classes
of equal size, by considering two embeddings of $G$ as equivalent if the
induced embeddings of $H$ are equal (up to rotations and translations). If
this assumption was violated, the resulting lower bound would be even better;
thus we can safely state the following bound:
\begin{equation}\label{eq:bound_genfan}
  \maxlamIIn \geq 2^{(n-|W|)\modop(|V|-|W|)} \cdot\lamII(H) \cdot
  \left(\frac{\lamII(G)}{\lamII(H)}\right)^{\!\lfloor(n-|W|)/(|V|-|W|)\rfloor} \qquad (n\geq|W|).
\end{equation}
Note that the previously described fan construction is a special instance of the
generalized fan, by taking as the subgraph $H$ a triangle with $\lamII(H)=2$.
Similarly, also the caterpillar construction can be seen as a special case, by
taking for $H$ a graph with $2$ vertices and Laman number~$1$.
To indicate the subgraph of a generalized fan construction we also write $H$-fan.
Using our encoding for graphs the usual fan would be denoted by $7$-fan.
The fan fixing the 4-vertex Laman graph is then denoted by $31$-fan.
Table~\ref{figure:fanbases} shows these bases.

\subsection{Rigid graphs with many embeddings}\label{sec:graphs-2d}

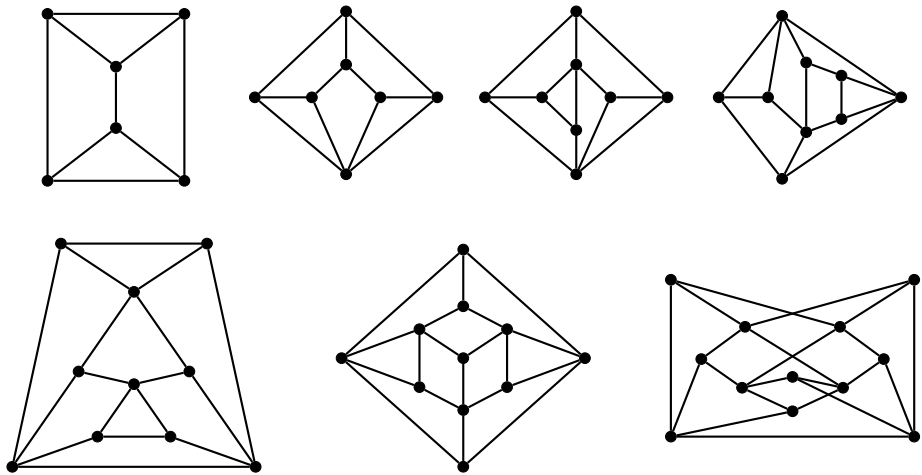
\begin{figure}[H]
	\begin{center}
		\begin{tabular}{c@{\hspace{12pt}}c@{\hspace{12pt}}c@{\hspace{12pt}}c}
			\begin{tikzpicture}[scale=0.06] 
				\draw[white] (-21,22) rectangle (22,-21);
				\vertex (a) at (-15.,-18.4482) {};
				\vertex (b) at (15,-18.4482) {};
				\vertex (c) at (0.,-6.78132) {};
				\vertex (d) at (-15.,18.4482) {};
				\vertex (e) at (15.,18.4482) {};
				\vertex (f) at (0.,6.78132) {};
				\draw[edge]  (a)edge(b) (a)edge(c) (a)edge(d) (b)edge(c) (b)edge(e) (c)edge(f) (d)edge(e) (d)edge(f) (e)edge(f);
			\end{tikzpicture}
			&
			\begin{tikzpicture}[scale=0.06] 
				\draw[white] (-21,22) rectangle (22,-21);
				\vertex (a) at (0.00,-17.00) {};
				\vertex (b) at (20.00,0) {};
				\vertex (c) at (7.50,0) {};
				\vertex (d) at (-20.00,0) {};
				\vertex (e) at (-7.50,0) {};
				\vertex (f) at (0.00,19.00) {};
				\vertex (g) at (0.00,7.24) {};
				\draw[edge]  (a)edge(b) (a)edge(c) (a)edge(d) (a)edge(e) (b)edge(c)
					(b)edge(f) (c)edge(g) (d)edge(e) (d)edge(f) (e)edge(g) (f)edge(g);
			\end{tikzpicture}
			&
			\begin{tikzpicture}[scale=0.06] 
				\draw[white] (-21,22) rectangle (22,-21);
				\vertex (a) at (-0.00,-17.00) {};
				\vertex (b) at (-0.00,7.24) {};
				\vertex (c) at (7.50,0.00) {};
				\vertex (d) at (-0.00,-7.24) {};
				\vertex (e) at (20.00,0.00) {};
				\vertex (f) at (-20.00,0.00) {};
				\vertex (g) at (-7.50,0.00) {};
				\vertex (h) at (-0.00,19.00) {};
				\draw[edge]  (a)edge(c) (a)edge(d) (a)edge(e) (a)edge(f) (b)edge(c)
					(b)edge(d) (b)edge(g) (b)edge(h) (c)edge(e) (d)edge(g) (e)edge(h)
					(f)edge(g) (f)edge(h);
			\end{tikzpicture}
			&
			\begin{tikzpicture}[scale=0.06] 
				\draw[white] (-21,22) rectangle (22,-21);
				\vertex (a) at (21.00,0.00) {};
				\vertex (b) at (-5.10,18.00) {};
				\vertex (c) at (0.17,-7.70) {};
				\vertex (d) at (7.90,-4.80) {};
				\vertex (e) at (-5.10,-18.00) {};
				\vertex (f) at (0.17,7.70) {};
				\vertex (g) at (-8.20,0.00) {};
				\vertex (h) at (7.90,4.80) {};
				\vertex (i) at (-19.00,0.00) {};
				\draw[edge]  (a)edge(b) (a)edge(d) (a)edge(e) (a)edge(h) (b)edge(f)
					(b)edge(g) (b)edge(i) (c)edge(d) (c)edge(e) (c)edge(f) (c)edge(g)
					(d)edge(h) (e)edge(i) (f)edge(h) (g)edge(i);
			\end{tikzpicture}
		\end{tabular}\\[12pt]
		\begin{tabular}{c@{\hspace{9pt}}c@{\hspace{9pt}}c}
			\begin{tikzpicture}[scale=0.08] 
				\draw[white] (-25,22) rectangle (25,-18);
				\vertex (a) at (-20.00,-16.00) {};
				\vertex (b) at (20.00,-16.00) {};
				\vertex (c) at (-0.00,-2.30) {};
				\vertex (d) at (-0.00,13.00) {};
				\vertex (e) at (-9.10,-0.20) {};
				\vertex (f) at (9.10,-0.20) {};
				\vertex (g) at (-6.00,-11.00) {};
				\vertex (h) at (-12.00,21.00) {};
				\vertex (i) at (6.00,-11.00) {};
				\vertex (j) at (12.00,21.00) {};
				\draw[edge]  (a)edge(b) (a)edge(e) (a)edge(g) (a)edge(h) (b)edge(f)
					(b)edge(i) (b)edge(j) (c)edge(e) (c)edge(f) (c)edge(g) (c)edge(i)
					(d)edge(e) (d)edge(f) (d)edge(h) (d)edge(j) (g)edge(i) (h)edge(j);
			\end{tikzpicture}
			&
			\begin{tikzpicture}[scale=0.08] 
				\draw[white] (-25,20) rectangle (25,-20);
				\vertex (a) at (7.20,4.80) {};
				\vertex (b) at (20.00,-0.00) {};
				\vertex (c) at (-7.20,4.80) {};
				\vertex (d) at (-20.00,-0.00) {};
				\vertex (e) at (0.00,-8.60) {};
				\vertex (f) at (7.20,-4.80) {};
				\vertex (g) at (0.00,-0.00) {};
				\vertex (h) at (0.00,-18.00) {};
				\vertex (i) at (-7.20,-4.80) {};
				\vertex (j) at (0.00,8.60) {};
				\vertex (k) at (0.00,18.00) {};
				\draw[edge]  (a)edge(b) (a)edge(f) (a)edge(g) (a)edge(j) (b)edge(f)
					(b)edge(h) (b)edge(k) (c)edge(d) (c)edge(g) (c)edge(i) (c)edge(j)
					(d)edge(h) (d)edge(i) (d)edge(k) (e)edge(f) (e)edge(g) (e)edge(h)
					(e)edge(i) (j)edge(k);
			\end{tikzpicture}
			&
			\begin{tikzpicture}[scale=0.08] 
				\draw[white] (-25,21) rectangle (25,-19);
				\vertex (a) at (20.00,-12.00) {};
				\vertex (b) at (-20.00,-12.00) {};
				\vertex (c) at (7.80,6.20) {};
				\vertex (d) at (-8.30,-3.90) {};
				\vertex (e) at (8.30,-3.90) {};
				\vertex (f) at (-7.80,6.20) {};
				\vertex (g) at (15.00,0.85) {};
				\vertex (h) at (20.00,14.00) {};
				\vertex (i) at (0.00,-2.10) {};
				\vertex (j) at (-20.00,14.00) {};
				\vertex (k) at (0.00,-7.80) {};
				\vertex (l) at (-15.00,0.85) {};
				\draw[edge] (a)edge(b) (a)edge(g) (a)edge(h) (a)edge(i) (b)edge(j)
				(b)edge(k) (b)edge(l) (c)edge(d) (c)edge(g) (c)edge(h) (c)edge(j)
				(d)edge(i) (d)edge(k) (d)edge(l) (e)edge(f) (e)edge(g) (e)edge(i)
				(e)edge(k) (f)edge(h) (f)edge(j) (f)edge(l);
			\end{tikzpicture}
		\end{tabular}
	\end{center}\vspace{-8pt}
	\caption{Unique Laman graphs with $6\leq n\leq12$ with maximal number of embeddings (see Table~\ref{table:max_laman_numbers}, encodings see Table~\ref{table:enc:max_laman_numbers}).}
	\label{figure:max_graphs}
\end{figure}

In order to get good lower bounds, we need particular Laman graphs
that have a large number of embeddings. For this purpose we have computed the Laman numbers
of all Laman graphs with up to $n=12$ vertices. 
We did so using the algorithm of \cite{SCgroup} (see \cite{SCgroupElectronicMaterial} for
an implementation and \cite{EurocombVersion} for a streamlined extended abstract).
For each $3 \leq n \leq 12$ we have
identified the (unique) Laman graph with the highest number of
embeddings. We present these numbers in
Table~\ref{table:max_laman_numbers} and the corresponding graphs for
$6\leq n\leq12$ appear in Figure~\ref{figure:max_graphs}.

\begin{table}[H]
	\begin{center}
		\begin{tabular}{lcccccccc}
			\toprule
			$n$
			& 6
			& 7
			& 8
			& 9
			& 10
			& 11
			& 12
			\\
			min
			& 16
			& 32
			& \phantom{1}64
			& 128
			& 256
			& \phantom{2}512
			& 1024
			\\
			$\maxlamIIn$
			& 24
			& 56
			& 136
			& 344
			& 880
			& 2288
			& 6180
			\\
			lower
			& 24
			& 48
			& \phantom{1}96
			& 288
			& 576
			& -
			& -
			\\
			\bottomrule
		\end{tabular}
	\end{center}
	\caption{Minimal and maximal Laman number among all $n$\hbox{-}vertex Laman graphs;
	the minimum is $2^{n-2}$ and it is achieved, for example, on Laman graphs that are
	constructible by using only Henneberg steps of type~1. The row labeled with ``lower'' contains
	the bounds from~\cite{Emiris2009}.}
	\label{table:max_laman_numbers}
\end{table}

There are 44\,176\,717 Laman graphs with 12 vertices, and therefore it was a
major undertaking to compute the Laman numbers of all of them; it took 56
processor days to complete this task. Hence it is unrealistic to do the same
for all Laman graphs with 13 or more vertices. In order to proceed further, we
developed some heuristics to construct graphs with very high Laman numbers,
albeit not necessarily the highest one.  The properties that we formulate for
the families $T(n)$ and $S(n)$ below are inspired by inspecting the few known
graphs that achieve the maximal Laman number~$\Mtwo(n)$ (see
Figure~\ref{figure:max_graphs} and Table~\ref{table:max_laman_numbers}).
More precisely, we consider the set $T(n)$ of Laman graphs with $n$ vertices
that satisfy the following additional properties.

\begin{definition}\label{def:special_graphs}
	We say that a Laman graph $G=(V,E)$ with $n$ vertices is an element of $T(n)$ iff
	\begin{itemize}
		\item $G$ is a planar graph, that is it can be embedded in the plane without
			crossings of edges.
		\item Each vertex of $G$ has degree 3 or 4; in this case the Laman condition
			(Definition~\ref{definition:laman_graph}) implies that there are
			exactly 6 vertices of degree 3 and $|V|-\text{6}$ vertices of degree 4.
		\item There are precisely two 3-cycles, and the number of 4 cycles is $|V|-\text{3}$. Note that 
			we count only nontrivial cycles all of whose edges are distinct.
			Moreover, the 3-cycles are disjoint, that is they do not share an edge.
			By Euler's formula the number of faces (including the outer, unlimited one)
			is given by $\text{2}-|V|+|E|=|V|-\text{1}$, and hence each of the cycles
			is the boundary of a face.
	\end{itemize}
\end{definition}

These properties are quite selective: for example, the set $T(12)$ contains
only $18$ (out of 44 million!) Laman graphs, and the set $T(18)$ has the
manageable cardinality~$188$. For $n\leq11$ we have that
$\max_{G\in T(n)}\bigl(\lamII(G)\bigr)=\Mtwo(n)$. In contrast, the $12$-vertex graph
with the highest Laman number is not in $T(12)$, since it is not
planar and does not have any $3$-cycles. Nevertheless, it satisfies
the condition on the vertex degrees.  Furthermore, the graph with the
highest Laman number in $T(12)$ is the graph with the second highest
Laman number with 12 vertices.  Hence, it is also the graph with the
highest Laman number which does contain a 3-cycle. For
$13\leq n\leq18$ we have constructed all Laman graphs $T(n)$ and among
them identified the one with the highest Laman number. We summarize
the results in Table~\ref{table:big_laman_numbers}; the
corresponding graphs are displayed in Figure~\ref{figure:max_graphs}.

\begin{table}[H]
	\begin{center}
		\begin{tabular}{lcccccccc}
			\toprule
			$n$
			& 12
			& 13
			& 14
			& 15
			& 16
			& 17
			& 18
			\\
			$M_T(n)$
			& 5952
			& 15056
			& 39696
			& 105384
			& 277864
			& 731336
			& 1953816
			\\
			$M_S(n)$
			& 6180
			& 15536
			& 42780
			& 112752
			& 312636
			& 870414
			& 2237312\\
			\bottomrule
		\end{tabular}
	\end{center}\vspace{-4pt}
	\caption{With $M_T(n)$ we denote the maximal Laman number of the graphs in~$T(n)$.
    In the row below we give the highest Laman numbers that we have found so far
    by looking at graphs in~$S(n)$ (exhaustive for $n\leq15$ but incomplete for $n>15$).}
	\label{table:big_laman_numbers}
\end{table}

\begin{figure}[ht]
  \begin{center}
    \begin{tabular}{c@{\hspace{12pt}}c@{\hspace{12pt}}c@{\hspace{12pt}}c}
      \begin{tikzpicture}[scale=0.08] 
        \draw[white] (-21,22) rectangle (22,-21);
        \vertex (a) at (0.00,-18.00) {};
        \vertex (b) at (20.00,0.53) {};
        \vertex (c) at (-20.00,0.53) {};
        \vertex (d) at (8.20,-4.20) {};
        \vertex (e) at (-2.70,-6.40) {};
        \vertex (f) at (0.00,8.20) {};
        \vertex (g) at (-8.20,-4.20) {};
        \vertex (h) at (2.70,-6.40) {};
        \vertex (i) at (0.00,19.00) {};
        \vertex (j) at (8.10,5.10) {};
        \vertex (k) at (-8.10,5.10) {};
        \vertex (l) at (0.00,0.53) {};
        \draw[edge] (a)edge(b) (a)edge(c) (a)edge(g) (a)edge(h) (b)edge(d)
        (b)edge(i) (b)edge(j) (c)edge(g) (c)edge(i) (c)edge(k) (d)edge(h)
        (d)edge(j) (d)edge(l) (e)edge(g) (e)edge(h) (e)edge(k) (e)edge(l)
        (f)edge(i) (f)edge(j) (f)edge(k) (f)edge(l);
      \end{tikzpicture}
      &
      \begin{tikzpicture}[scale=0.08] 
	\draw[white] (-21,22) rectangle (22,-21);
	\vertex (a) at (20.00,17.00) {};
	\vertex (b) at (20.00,-16.00) {};
	\vertex (c) at (8.00,10.00) {};
	\vertex (d) at (-20.00,17.00) {};
	\vertex (e) at (8.00,-9.00) {};
	\vertex (f) at (-8.00,-9.00) {};
	\vertex (g) at (-4.20,0.48) {};
	\vertex (h) at (15.00,0.48) {};
	\vertex (i) at (-20.00,-16.00) {};
	\vertex (j) at (-8.00,10.00) {};
	\vertex (k) at (4.20,0.48) {};
	\vertex (l) at (-15.00,0.48) {};
	\vertex (m) at (0.00,-5.70) {};
	\draw[edge] (a)edge(b) (a)edge(c) (a)edge(d) (a)edge(h) (b)edge(e)
	(b)edge(f) (b)edge(i) (c)edge(h) (c)edge(j) (c)edge(k) (d)edge(i)
	(d)edge(j) (d)edge(l) (e)edge(h) (e)edge(k) (e)edge(m) (f)edge(i)
	(f)edge(l) (f)edge(m) (g)edge(j) (g)edge(k) (g)edge(l) (g)edge(m);
      \end{tikzpicture}
      &
      \begin{tikzpicture}[scale=0.08] 
	\draw[white] (-21,22) rectangle (22,-21);
	\vertex (a) at (-19.00,0.00) {};
	\vertex (b) at (0.64,-3.70) {};
	\vertex (c) at (0.64,-18.00) {};
	\vertex (d) at (-4.20,-7.70) {};
	\vertex (e) at (-8.30,3.60) {};
	\vertex (f) at (0.64,18.00) {};
	\vertex (g) at (12.00,0.00) {};
	\vertex (h) at (0.64,3.70) {};
	\vertex (i) at (-8.30,-3.60) {};
	\vertex (j) at (5.40,-7.70) {};
	\vertex (k) at (21.00,0.00) {};
	\vertex (l) at (-2.10,0.00) {};
	\vertex (m) at (-4.20,7.70) {};
	\vertex (n) at (5.40,7.70) {};
	\draw[edge] (a)edge(c) (a)edge(e) (a)edge(f) (a)edge(i) (b)edge(d)
	(b)edge(g) (b)edge(h) (b)edge(j) (c)edge(d) (c)edge(j) (c)edge(k)
	(d)edge(i) (d)edge(l) (e)edge(i) (e)edge(l) (e)edge(m) (f)edge(k)
	(f)edge(m) (f)edge(n) (g)edge(j) (g)edge(k) (g)edge(n) (h)edge(l)
	(h)edge(m) (h)edge(n);
      \end{tikzpicture}
      &
      \begin{tikzpicture}[scale=0.08] 
	\draw[white] (-21,22) rectangle (22,-21);
	\vertex (a) at (-20.00,-18.00) {};
	\vertex (b) at (0.37,13.00) {};
	\vertex (c) at (-20.00,18.00) {};
	\vertex (d) at (20.00,-18.00) {};
	\vertex (e) at (-5.30,4.10) {};
	\vertex (f) at (0.31,-12.00) {};
	\vertex (g) at (6.10,-3.80) {};
	\vertex (h) at (-15.00,-3.50) {};
	\vertex (i) at (13.00,9.20) {};
	\vertex (j) at (-4.60,-12.00) {};
	\vertex (k) at (6.10,3.80) {};
	\vertex (l) at (20.00,18.00) {};
	\vertex (m) at (-13.00,9.20) {};
	\vertex (n) at (15.00,-3.50) {};
	\vertex (o) at (-5.50,-3.80) {};
	\draw[edge] (a)edge(c) (a)edge(d) (a)edge(h) (a)edge(j) (b)edge(c)
	(b)edge(e) (b)edge(i) (b)edge(k) (c)edge(l) (c)edge(m) (d)edge(f)
	(d)edge(l) (d)edge(n) (e)edge(g) (e)edge(m) (e)edge(o) (f)edge(g)
	(f)edge(j) (f)edge(o) (g)edge(k) (g)edge(n) (h)edge(j) (h)edge(m)
	(h)edge(o) (i)edge(k) (i)edge(l) (i)edge(n);
      \end{tikzpicture}
    \end{tabular}
    \begin{tabular}{c@{\hspace{12pt}}c@{\hspace{12pt}}c}
      \begin{tikzpicture}[scale=0.08] 
	\draw[white] (-21,22) rectangle (22,-21);
	\vertex (a) at (5.00,-12.00) {};
	\vertex (b) at (-15.00,16.00) {};
	\vertex (c) at (20.00,-16.00) {};
	\vertex (d) at (6.00,-6.10) {};
	\vertex (e) at (0.00,-1.90) {};
	\vertex (f) at (-20.00,-16.00) {};
	\vertex (g) at (0.00,10.00) {};
	\vertex (h) at (-13.00,3.10) {};
	\vertex (i) at (13.00,3.10) {};
	\vertex (j) at (-6.00,-6.10) {};
	\vertex (k) at (-5.00,-12.00) {};
	\vertex (l) at (15.00,16.00) {};
	\vertex (m) at (13.00,-6.90) {};
	\vertex (n) at (6.00,2.30) {};
	\vertex (o) at (-6.00,2.30) {};
	\vertex (p) at (-13.00,-6.90) {};
	\draw[edge] (a)edge(c) (a)edge(d) (a)edge(e) (a)edge(k) (b)edge(f)
	(b)edge(g) (b)edge(h) (b)edge(l) (c)edge(f) (c)edge(l) (c)edge(m)
	(d)edge(i) (d)edge(m) (d)edge(n) (e)edge(j) (e)edge(n) (e)edge(o)
	(f)edge(k) (f)edge(p) (g)edge(i) (g)edge(n) (g)edge(o) (h)edge(j)
	(h)edge(o) (h)edge(p) (i)edge(l) (i)edge(m) (j)edge(k) (j)edge(p);
      \end{tikzpicture}
      &
      \begin{tikzpicture}[scale=0.08] 
	\draw[white] (-21,22) rectangle (22,-21);
	\vertex (a) at (-20.00,-16.00) {};
	\vertex (b) at (-20.00,16.00) {};
	\vertex (c) at (11.00,-9.00) {};
	\vertex (d) at (-6.40,4.60) {};
	\vertex (e) at (-1.70,-12.00) {};
	\vertex (f) at (-15.00,-0.34) {};
	\vertex (g) at (0.08,11.00) {};
	\vertex (h) at (20.00,16.00) {};
	\vertex (i) at (4.80,-1.90) {};
	\vertex (j) at (15.00,-0.34) {};
	\vertex (k) at (-4.60,-1.90) {};
	\vertex (l) at (20.00,-16.00) {};
	\vertex (m) at (-11.00,7.60) {};
	\vertex (n) at (-11.00,-9.00) {};
	\vertex (o) at (0.08,-6.10) {};
	\vertex (p) at (12.00,11.00) {};
	\vertex (q) at (6.60,4.50) {};
	\draw[edge] (a)edge(b) (a)edge(e) (a)edge(f) (a)edge(l) (b)edge(g)
	(b)edge(h) (b)edge(m) (c)edge(e) (c)edge(i) (c)edge(j) (c)edge(l)
	(d)edge(g) (d)edge(i) (d)edge(k) (d)edge(m) (e)edge(n) (e)edge(o)
	(f)edge(k) (f)edge(m) (f)edge(n) (g)edge(p) (g)edge(q) (h)edge(j)
	(h)edge(l) (h)edge(p) (i)edge(o) (i)edge(q) (j)edge(p) (j)edge(q)
	(k)edge(n) (k)edge(o);
      \end{tikzpicture}
      &
      \begin{tikzpicture}[scale=0.08] 
	\draw[white] (-21,22) rectangle (22,-21);
	\vertex (a) at (-20.00,-17.00) {};
	\vertex (b) at (-11.00,-11.00) {};
	\vertex (c) at (-16.00,0.00) {};
	\vertex (d) at (15.00,4.80) {};
	\vertex (e) at (7.90,0.00) {};
	\vertex (f) at (7.20,8.50) {};
	\vertex (g) at (20.00,-17.00) {};
	\vertex (h) at (-3.20,-4.80) {};
	\vertex (i) at (-11.00,11.00) {};
	\vertex (j) at (20.00,17.00) {};
	\vertex (k) at (7.20,-8.50) {};
	\vertex (l) at (-3.20,4.80) {};
	\vertex (m) at (-20.00,17.00) {};
	\vertex (n) at (0.22,-13.00) {};
	\vertex (o) at (-7.70,0.00) {};
	\vertex (p) at (15.00,-4.80) {};
	\vertex (q) at (0.10,0.00) {};
	\vertex (r) at (0.22,13.00) {};
	\draw[edge] (a)edge(b) (a)edge(c) (a)edge(g) (a)edge(m) (b)edge(c)
	(b)edge(h) (b)edge(n) (c)edge(i) (c)edge(o) (d)edge(e) (d)edge(f)
	(d)edge(j) (d)edge(p) (e)edge(f) (e)edge(k) (e)edge(q) (f)edge(l)
	(f)edge(r) (g)edge(j) (g)edge(n) (g)edge(p) (h)edge(k) (h)edge(o)
	(h)edge(q) (i)edge(l) (i)edge(m) (i)edge(r) (j)edge(m) (j)edge(r)
	(k)edge(n) (k)edge(p) (l)edge(o) (l)edge(q);
      \end{tikzpicture}
    \end{tabular}
  \end{center}
  \caption{Laman graphs in $T(n)$ with $12\leq n\leq18$ vertices; for each $n$ the
    graph with the largest Laman number among the Laman graphs in $T(n)$ is displayed. The corresponding Laman
    numbers are given in Table~\ref{table:big_laman_numbers} (encodings see Table~\ref{table:enc:max_graphs_triangle}).}
  \label{figure:max_graphs_triangle}
\end{figure}
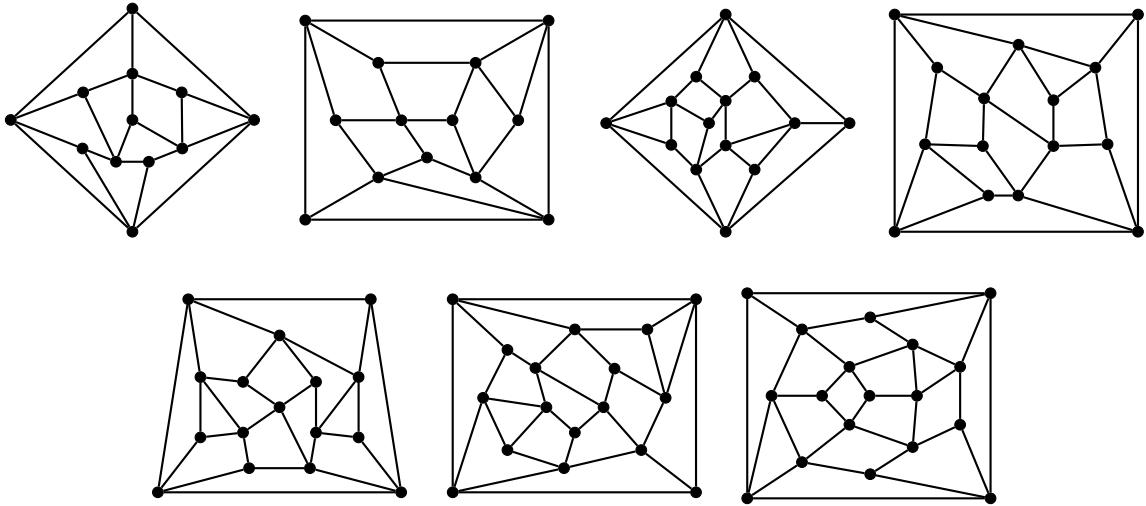

\begin{figure}[H]
  \begin{center}
		\begin{tabular}{c@{\qquad}c@{\qquad}c@{\qquad}c}
			\begin{tikzpicture}
				\node[vertex] (1) at (-0.5, -0.866025) {};
				\node[vertex] (2) at (-0.5, 0.866025) {};
				\node[vertex] (3) at (0.5, 0.866025) {};
				\node[vertex] (4) at (1., 0.) {};
				\node[vertex] (5) at (0.5, -0.866025) {};
				\node[vertex] (6) at (-1., 0.) {};
				\draw[edge] (1)edge(4) (1)edge(5) (1)edge(6) (2)edge(3) (2)edge(5)
				(2)edge(6) (3)edge(4) (3)edge(6) (4)edge(5) ;,
			\end{tikzpicture}
			&
			\begin{tikzpicture}
				\node[vertex] (1)	at (-0.433884, -0.900969) {};
				\node[vertex] (2) at (0.433884, -0.900969) {};
				\node[vertex] (3) at (0.974928, -0.222521) {};
				\node[vertex] (4) at (-0.781831, 0.62349) {};
				\node[vertex] (5) at (-0.974928, -0.222521) {};
				\node[vertex] (6) at (0.781831, 0.62349) {};
				\node[vertex] (7) at (0., 1.) {};
				\draw[edge] (1)edge(2) (1)edge(5) (1)edge(6) (2)edge(3) (2)edge(4)
				(3)edge(6) (3)edge(7) (4)edge(5) (4)edge(7) (5)edge(7) (6)edge(7) ;
			\end{tikzpicture}
			&
			\begin{tikzpicture}
				\node[vertex] (1) at (0.382683, -0.92388) {};
				\node[vertex] (2) at (-0.382683, 0.92388) {};
				\node[vertex] (3) at (-0.92388, 0.382683) {};
				\node[vertex] (4) at (0.92388, -0.382683) {};
				\node[vertex] (5) at (0.92388, 0.382683) {};
				\node[vertex] (6) at (-0.92388, -0.382683) {};
				\node[vertex] (7) at (0.382683, 0.92388) {};
				\node[vertex] (8) at (-0.382683, -0.92388) {};
				\draw[edge] (1)edge(2) (1)edge(4) (1)edge(8) (2)edge(3) (2)edge(7)
				(3)edge(6) (3)edge(8) (4)edge(5) (4)edge(7) (5)edge(7) (5)edge(8)
				(6)edge(7) (6)edge(8);
			\end{tikzpicture}
			&
			\begin{tikzpicture}
				\node[vertex] (1) at (-0.984808,0.173648) {};
				\node[vertex] (2) at (0.984808, 0.173648) {};
				\node[vertex] (3) at (0.642788, 0.766044) {};
				\node[vertex] (4) at (-0.866025, -0.5) {};
				\node[vertex] (5) at (0.866025, -0.5) {};
				\node[vertex] (6) at (-0.642788, 0.766044) {};
				\node[vertex] (7) at (0., 1.) {};
				\node[vertex] (8) at (0.34202, -0.939693) {};
				\node[vertex] (9) at (-0.34202, -0.939693) {};
				\draw[edge] (1)edge(4) (1)edge(6) (1)edge(9) (2)edge(3) (2)edge(5)
				(2)edge(8) (3)edge(7) (3)edge(9) (4)edge(7) (4)edge(9) (5)edge(7)
				(5)edge(8) (6)edge(7) (6)edge(8) (8)edge(9);
			\end{tikzpicture}
		\end{tabular}
		\vspace{16pt}\\
    \begin{tabular}{c@{\qquad}c@{\qquad}c}
			\begin{tikzpicture}
				\node[vertex] (1)	at (-0.809017, 0.587785) {};
				\node[vertex] (2) at (1., 0.) {};
				\node[vertex] (3) at (0.809017, -0.587785) {};
				\node[vertex] (4) at (-1., 0.) {};
				\node[vertex] (5) at (0.309017, 0.951057) {};
				\node[vertex] (6) at (-0.309017, -0.951057) {};
				\node[vertex] (7) at (0.809017, 0.587785) {};
				\node[vertex] (8) at (-0.809017, -0.587785) {};
				\node[vertex] (9) at (0.309017, -0.951057) {};
				\node[vertex] (10) at (-0.309017, 0.951057) {};
				\draw[edge] (1)edge(4) (1)edge(8) (1)edge(10) (2)edge(3) (2)edge(7)
				(2)edge(10) (3)edge(7) (3)edge(9) (4)edge(8) (4)edge(9) (5)edge(7)
				(5)edge(8) (5)edge(10) (6)edge(7) (6)edge(8) (6)edge(9) (9)edge(10);
			\end{tikzpicture}
			&
			\begin{tikzpicture}
				\node[vertex] (1) at (-0.281733, -0.959493) {};
				\node[vertex] (2) at (0.281733, -0.959493) {};
				\node[vertex] (3) at (-0.989821, -0.142315) {};
				\node[vertex] (4) at (0.540641, 0.841254) {};
				\node[vertex] (5) at (-0.540641, 0.841254) {};
				\node[vertex] (6) at (0.989821, -0.142315) {};
				\node[vertex] (7) at (0., 1.) {};
				\node[vertex] (8) at (0.75575, -0.654861) {};
				\node[vertex] (9) at (-0.909632, 0.415415) {};
				\node[vertex] (10) at (-0.75575, -0.654861) {};
				\node[vertex] (11) at (0.909632, 0.415415) {};
				\draw[edge] (1)edge(2) (1)edge(10) (1)edge(11) (2)edge(8) (2)edge(9)
				(3)edge(7) (3)edge(9) (3)edge(10) (4)edge(7) (4)edge(10) (4)edge(11)
				(5)edge(7) (5)edge(8) (5)edge(9) (6)edge(7) (6)edge(8) (6)edge(11)
				(8)edge(11) (9)edge(10);
			\end{tikzpicture}
			&
			\begin{tikzpicture}
				\node[vertex] (1) at (0.707107,-0.707107) {};
				\node[vertex] (2) at (-0.707107, 0.707107) {};
				\node[vertex] (3) at (-0.965926, -0.258819) {};
				\node[vertex] (4) at (0.965926, 0.258819) {};
				\node[vertex] (5) at (-0.258819, -0.965926) {};
				\node[vertex] (6) at (0.258819, 0.965926) {};
				\node[vertex] (7) at (0.707107, 0.707107) {};
				\node[vertex] (8) at (-0.707107, -0.707107) {};
				\node[vertex] (9) at (-0.258819, 0.965926) {};
				\node[vertex] (10) at (0.258819, -0.965926) {};
				\node[vertex] (11) at (0.965926, -0.258819) {};
				\node[vertex] (12) at (-0.965926, 0.258819) {};
				\draw[edge] (1)edge(10) (1)edge(11) (1)edge(12) (2)edge(9) (2)edge(11)
				(2)edge(12) (3)edge(8) (3)edge(10) (3)edge(12) (4)edge(7) (4)edge(9)
				(4)edge(11) (5)edge(7) (5)edge(8) (5)edge(10) (6)edge(7) (6)edge(8)
				(6)edge(9) (7)edge(12) (8)edge(11) (9)edge(10) ;
			\end{tikzpicture}
		\end{tabular}
  \end{center}
  \caption{Circular embedding of the Laman graphs with maximal Laman numbers for $6\leq n\leq12$.
    Note that these are the same graphs that are displayed in Figure~\ref{figure:max_graphs}.}
  \label{fig:symmax}
\end{figure}
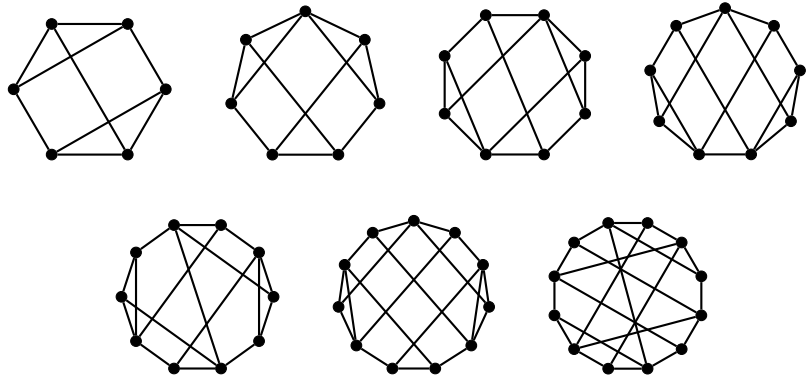

We have seen that for 12 vertices the maximal graph in $T(12)$ is not the one with the highest Laman number.
The same holds true for $13\leq n\leq 18$, which can be seen by looking at a different family of graphs:
We observed that the graphs which are known to be maximal according to their Laman number are Hamiltonian, 
i.e., they contain a path that visits each vertex exactly once (Hamiltonian path).
Hence, we focus on Hamiltonian graphs. The problem is that they cover still around $2/3$ of all Laman graphs (at least for small~$n$).
Therefore, we considered other properties of the known graphs with maximal Laman number.
One of these properties is the symmetry of a certain embedding.
\begin{definition}
  We say that a Laman graph $G=(V,E)$ is an element of $S(n)$ iff
  \begin{itemize}
    \item $G$ is Hamiltonian, i.e.\ it contains a Hamiltonian cycle $H$.
    \item There exists a circular embedding, i.e.\ an embedding $\rho$ such that $\rho(v)$ lies on the same circle for all $v\in V$, and $H$ is embedded on a regular $n$-gon.
    \item The figure obtained by the embedding is point resp.\ line symmetric for an even resp.\ odd number of vertices.
  \end{itemize}
\end{definition}
One can see that the maximal Laman graphs up to 12 vertices fulfill these symmetry properties (Figure~\ref{fig:symmax}).

\begin{figure}[H]
	\begin{center}
		\begin{tabular}{cccccc}
			\begin{tikzpicture} 
				\node[vertex] (1) at (-0.663123, -0.748511) {};
				\node[vertex] (2) at (-0.992709, 0.120537) {};
				\node[vertex] (3) at (0.992709, 0.120537) {};
				\node[vertex] (4) at (0.663123, -0.748511) {};
				\node[vertex] (5) at (-0.464723, 0.885456) {};
				\node[vertex] (6) at (0.464723, 0.885456) {};
				\node[vertex] (7) at (0.935016, -0.354605) {};
				\node[vertex] (8) at (-0.935016, -0.354605) {};
				\node[vertex] (9) at (0., 1.) {};
				\node[vertex] (10) at (0.239316, -0.970942) {};
				\node[vertex] (11) at (-0.239316, -0.970942) {};
				\node[vertex] (12) at (0.822984, 0.568065) {};
				\node[vertex] (13) at (-0.822984, 0.568065) {};
				\draw[edge] (1)edge(8) (1)edge(11) (1)edge(13) (2)edge(8) (2)edge(10)
				(2)edge(13) (3)edge(7) (3)edge(11) (3)edge(12) (4)edge(7) (4)edge(10)
				(4)edge(12) (5)edge(9) (5)edge(10) (5)edge(13) (6)edge(9) (6)edge(11)
				(6)edge(12) (7)edge(9) (7)edge(13) (8)edge(9) (8)edge(12)
				(10)edge(11) ;
			\end{tikzpicture}
			&
			\begin{tikzpicture} 
				\node[vertex] (1) at (1., 0.) {};
				\node[vertex] (2) at (-0.62349, 0.781831) {};
				\node[vertex] (3) at (0.222521, 0.974928) {};
				\node[vertex] (4) at (-1., 0.) {};
				\node[vertex] (5) at (-0.222521, -0.974928) {};
				\node[vertex] (6) at (0.62349, -0.781831) {};
				\node[vertex] (7) at (-0.62349, -0.781831) {};
				\node[vertex] (8) at (0.62349, 0.781831) {};
				\node[vertex] (9) at (-0.900969, 0.433884) {};
				\node[vertex] (10) at (-0.900969, -0.433884) {};
				\node[vertex] (11) at (0.900969, 0.433884) {};
				\node[vertex] (12) at (0.900969, -0.433884) {};
				\node[vertex] (13) at (0.222521, -0.974928) {};
				\node[vertex] (14) at (-0.222521, 0.974928) {};
				\draw[edge] (1)edge(11) (1)edge(12) (1)edge(14) (2)edge(9) (2)edge(10)
				(2)edge(14) (3)edge(7) (3)edge(8) (3)edge(14) (4)edge(9) (4)edge(10)
				(4)edge(13) (5)edge(7) (5)edge(8) (5)edge(13) (6)edge(11) (6)edge(12)
				(6)edge(13) (7)edge(10) (7)edge(12) (8)edge(9) (8)edge(11)
				(9)edge(12) (10)edge(11) (13)edge(14) ;
			\end{tikzpicture}
			&
			\begin{tikzpicture} 
				\node[vertex] (1) at (0.587785, -0.809017) {};
				\node[vertex] (2) at (-0.587785, -0.809017) {};
				\node[vertex] (3) at (0.406737, 0.913545) {};
				\node[vertex] (4) at (-0.406737, 0.913545) {};
				\node[vertex] (5) at (0.994522, -0.104528) {};
				\node[vertex] (6) at (-0.994522, -0.104528) {};
				\node[vertex] (7) at (-0.951057, 0.309017) {};
				\node[vertex] (8) at (0.951057, 0.309017) {};
				\node[vertex] (9) at (0., 1.) {};
				\node[vertex] (10) at (0.743145, 0.669131) {};
				\node[vertex] (11) at (-0.743145, 0.669131) {};
				\node[vertex] (12) at (-0.866025, -0.5) {};
				\node[vertex] (13) at (0.866025, -0.5) {};
				\node[vertex] (14) at (0.207912, -0.978148) {};
				\node[vertex] (15) at (-0.207912, -0.978148) {};
				\draw[edge] (1)edge(10) (1)edge(13) (1)edge(14) (2)edge(11) (2)edge(12)
				(2)edge(15) (3)edge(9) (3)edge(10) (3)edge(15) (4)edge(9) (4)edge(11)
				(4)edge(14) (5)edge(8) (5)edge(13) (5)edge(15) (6)edge(7) (6)edge(12)
				(6)edge(14) (7)edge(8) (7)edge(11) (7)edge(13) (8)edge(10)
				(8)edge(12) (9)edge(12) (9)edge(13) (10)edge(11) (14)edge(15) ;
			\end{tikzpicture}
			&
			\begin{tikzpicture} 
				\node[vertex] (1) at (-0.980785, -0.19509) {};
				\node[vertex] (2) at (-0.19509, -0.980785) {};
				\node[vertex] (3) at (0.19509, 0.980785) {};
				\node[vertex] (4) at (0.980785, 0.19509) {};
				\node[vertex] (5) at (-0.55557, 0.83147) {};
				\node[vertex] (6) at (0.55557, -0.83147) {};
				\node[vertex] (7) at (0.83147, -0.55557) {};
				\node[vertex] (8) at (-0.83147, 0.55557) {};
				\node[vertex] (9) at (0.83147, 0.55557) {};
				\node[vertex] (10) at (-0.19509, 0.980785) {};
				\node[vertex] (11) at (0.19509, -0.980785) {};
				\node[vertex] (12) at (-0.83147, -0.55557) {};
				\node[vertex] (13) at (-0.55557, -0.83147) {};
				\node[vertex] (14) at (-0.980785, 0.19509) {};
				\node[vertex] (15) at (0.980785, -0.19509) {};
				\node[vertex] (16) at (0.55557, 0.83147) {};
				\draw[edge] (1)edge(12) (1)edge(14) (1)edge(15) (2)edge(11) (2)edge(13)
				(2)edge(16) (3)edge(10) (3)edge(13) (3)edge(16) (4)edge(9)
				(4)edge(14) (4)edge(15) (5)edge(8) (5)edge(10) (5)edge(12) (6)edge(7)
				(6)edge(9) (6)edge(11) (7)edge(8) (7)edge(13) (7)edge(15) (8)edge(14)
				(8)edge(16) (9)edge(12) (9)edge(16) (10)edge(11) (10)edge(15)
				(11)edge(14) (12)edge(13) ;
			\end{tikzpicture}
			&
			\begin{tikzpicture} 
				\node[vertex] (1) at (0., 1.) {};
				\node[vertex] (2) at (0.361242, 0.932472) {};
				\node[vertex] (3) at (0.673696, 0.739009) {};
				\node[vertex] (4) at (0.895163, 0.445738) {};
				\node[vertex] (5) at (0.995734, 0.0922684) {};
				\node[vertex] (6) at (0.961826, -0.273663) {};
				\node[vertex] (7) at (0.798017, -0.602635) {};
				\node[vertex] (8) at (0.526432, -0.850217) {};
				\node[vertex] (9) at (0.18375, -0.982973) {};
				\node[vertex] (10) at (-0.18375, -0.982973) {};
				\node[vertex] (11) at (-0.526432, -0.850217) {};
				\node[vertex] (12) at (-0.798017, -0.602635) {};
				\node[vertex] (13) at (-0.961826, -0.273663) {};
				\node[vertex] (14) at (-0.995734, 0.0922684) {};
				\node[vertex] (15) at (-0.895163, 0.445738) {};
				\node[vertex] (16) at (-0.673696, 0.739009) {};
				\node[vertex] (17) at (-0.361242, 0.932472) {};
				\draw[edge] (13)edge(12) (13)edge(14) (13)edge(1) (6)edge(7) (6)edge(5)
				(6)edge(1) (15)edge(2) (15)edge(14) (15)edge(16) (4)edge(17)
				(4)edge(5) (4)edge(3) (8)edge(9) (8)edge(7) (8)edge(16) (11)edge(10)
				(11)edge(12) (11)edge(3) (2)edge(9) (2)edge(3) (2)edge(1)
				(17)edge(10) (17)edge(16) (17)edge(1) (10)edge(9) (10)edge(14)
				(9)edge(5) (12)edge(5) (12)edge(16) (7)edge(14) (7)edge(3) ;
			\end{tikzpicture}
			&
			\begin{tikzpicture} 
			  \node[vertex] (1) at (0.5, 0.866025) {};
				\node[vertex] (2) at (-0.5, -0.866025) {};
				\node[vertex] (3) at (-1., 0.) {};
				\node[vertex] (4) at (1., 0.) {};
				\node[vertex] (5) at (-0.5, 0.866025) {};
				\node[vertex] (6) at (0.5, -0.866025) {};
				\node[vertex] (7) at (0.766044, -0.642788) {};
				\node[vertex] (8) at (-0.766044, 0.642788) {};
				\node[vertex] (9) at (0.173648, -0.984808) {};
				\node[vertex] (10) at (-0.173648, 0.984808) {};
				\node[vertex] (11) at (0.939693, 0.34202) {};
				\node[vertex] (12) at (-0.939693, -0.34202) {};
				\node[vertex] (13) at (-0.939693, 0.34202) {};
				\node[vertex] (14) at (0.939693, -0.34202) {};
				\node[vertex] (15) at (0.173648, 0.984808) {};
				\node[vertex] (16) at (-0.173648, -0.984808) {};
				\node[vertex] (17) at (-0.766044, -0.642788) {};
				\node[vertex] (18) at (0.766044, 0.642788) {};
				\draw[edge] (1)edge(15) (1)edge(16) (1)edge(18) (2)edge(15) (2)edge(16) 
				(2)edge(17) (3)edge(12) (3)edge(13) (3)edge(18) (4)edge(11) 
				(4)edge(14) (4)edge(17) (5)edge(8) (5)edge(10) (5)edge(14) (6)edge(7) 
				(6)edge(9) (6)edge(13) (7)edge(12) (7)edge(14) (7)edge(18) 
				(8)edge(11) (8)edge(13) (8)edge(17) (9)edge(10) (9)edge(11)
				(9)edge(16) (10)edge(12) (10)edge(15) (11)edge(18) (12)edge(17)
				(13)edge(15) (14)edge(16);
			\end{tikzpicture}
		\end{tabular}
	\end{center}
	\caption{For $n=13,\dots,18$ we display the graph from $S(n)$ with the highest Laman number
          (given in Table~\ref{table:big_laman_numbers}) found so far (encodings see Table~\ref{table:enc:max_graphs_high}).}
	\label{figure:max_graphs_high}
\end{figure}
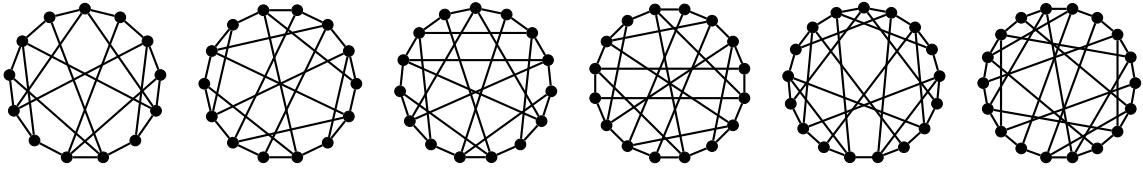

We computed the Laman numbers of all graphs in $S(n)$ up to $n=15$.
Unfortunately, for larger $n$ the set $S(n)$ still contains too many graphs.
For $n=15$ there are already $85058$ such graphs.
Performing the computations on a subset of~$S(n)$ yields the graphs shown in
Figure~\ref{figure:max_graphs_high}.

\subsection{Lower bounds}\label{sec:bounds-2d}

We now use these results to derive new and better lower bounds than the
previously known ones. We apply the caterpillar construction to the Laman
graphs with the maximal number of embeddings for $6\leq n\leq12$, and for
$13\leq n\leq18$ we use the graphs found by exploring the set $S(n)$ (see
Figure~\ref{figure:max_graphs_high} and Table~\ref{table:big_laman_numbers}).
The fan construction is applied to
the maximal Laman graphs for $6\leq n\leq11$ only, since it is not applicable
to the maximal graph with $12$ vertices (Figure~\ref{figure:max_graphs}).  Hence,
for the remaining cases, $12\leq n\leq18$, the fan construction is applied to
the maximal graph in $T(n)$.
In Table~\ref{table:bounds} the results obtained by these graphs are written in a separate column.
The results in the next column are obtained by randomly found graphs which contain a triangle and have a higher Laman number than the one in~$T(n)$.

For $7\leq n\leq11$ we also tried the generalized
fan construction: among all Laman graphs whose vertex degrees are at
least $3$---we can exclude Laman graphs that have vertices of degree $2$ since
they can be derived from a smaller graph by Henneberg steps of type~1, thereby only
doubling the embedding number---we selected all graphs that have the
$4$-vertex Laman graph as a subgraph. Then we computed their Laman numbers in
order to find the maximum that can be achieved among those graphs. Until 12 vertices
the lower bounds, according to~\eqref{eq:bound_genfan}, are not as good as
those obtained by the standard fan construction.
For higher $n$, randomly found graphs show improvements over the fan construction for the graphs we have found.
Note that since the graphs are found only randomly this does not show any results on whether the factors are indeed better.

\begin{figure}[H]
	\begin{center}
		\begin{tabular}{ccc}
			\begin{tikzpicture}[scale=0.05]
				\vertex (a) at (0.00,6.90) {};
				\vertex (b) at (0.00,22.00) {};
				\vertex (c) at (20.00,6.90) {};
				\vertex (d) at (-20.00,6.90) {};
				\vertex (e) at (0.00,-8.50) {};
				\vertex (f) at (14.00,-17.00) {};
				\vertex (g) at (-14.00,-17.00) {};
				\draw[edge] (a)edge(b) (a)edge(c) (a)edge(d) (a)edge(e) (b)edge(c)
					(b)edge(d) (c)edge(f) (d)edge(g) (e)edge(f) (e)edge(g) (f)edge(g);
			\end{tikzpicture}
			&
			\begin{tikzpicture}[scale=0.05]
				\vertex (a) at (0.00,4.70) {};
				\vertex (b) at (0.00,23.00) {};
				\vertex (c) at (20.00,9.20) {};
				\vertex (d) at (-20.00,9.20) {};
				\vertex (e) at (7.80,-7.30) {};
				\vertex (f) at (-7.80,-7.30) {};
				\vertex (g) at (18.00,-16.00) {};
				\vertex (h) at (-18.00,-16.00) {};
				\draw[edge] (a)edge(b) (a)edge(c) (a)edge(d) (a)edge(e) (a)edge(f)
					(b)edge(c) (b)edge(d) (c)edge(g) (d)edge(h) (e)edge(f) (e)edge(g)
					(f)edge(h) (g)edge(h);
			\end{tikzpicture}
			&
			\begin{tikzpicture}[scale=0.06]
				\vertex (a) at (-6.90,-4.40) {};
				\vertex (b) at (10.00,-15.00) {};
				\vertex (c) at (0.00,19.00) {};
				\vertex (d) at (6.90,-4.40) {};
				\vertex (e) at (-10.00,-15.00) {};
				\vertex (f) at (6.90,5.30) {};
				\vertex (g) at (20.00,5.20) {};
				\vertex (h) at (-6.90,5.30) {};
				\vertex (i) at (-20.00,5.20) {};
				\draw[edge] (a)edge(b) (a)edge(d) (a)edge(e) (a)edge(f) (b)edge(d)
					(b)edge(e) (b)edge(g) (c)edge(f) (c)edge(g) (c)edge(h) (c)edge(i)
					(d)edge(h) (e)edge(i) (f)edge(g) (h)edge(i);
			\end{tikzpicture}
			\\
			48 & 112 & 288\\
			\begin{tikzpicture}[scale=0.06]
				\vertex (a) at (11.00,0.40) {};
				\vertex (b) at (16.00,-13.00) {};
				\vertex (c) at (0.03,4.00) {};
				\vertex (d) at (-16.00,-13.00) {};
				\vertex (e) at (3.80,-6.60) {};
				\vertex (f) at (20.00,9.00) {};
				\vertex (g) at (-3.70,-6.60) {};
				\vertex (h) at (-11.00,0.04) {};
				\vertex (i) at (0.03,16.00) {};
				\vertex (j) at (-20.00,9.00) {};
				\draw[edge] (a)edge(b) (a)edge(c) (a)edge(e) (a)edge(f) (b)edge(d)
					(b)edge(e) (b)edge(f) (c)edge(g) (c)edge(h) (c)edge(i) (d)edge(g)
					(d)edge(h) (d)edge(j) (e)edge(g) (f)edge(i) (h)edge(j) (i)edge(j);
			\end{tikzpicture}
			&
			\begin{tikzpicture}[scale=0.06]
				\vertex (a) at (0,9.40) {};
				\vertex (b) at (-15,-16.00) {};
				\vertex (c) at (15,-16.00) {};
				\vertex (d) at (0,-2.30) {};
				\vertex (e) at (-9.1,0.04) {};
				\vertex (f) at (9.1,0.04) {};
				\vertex (g) at (-15,11.00) {};
				\vertex (h) at (15,11.00) {};
				\vertex (i) at (-4.1,-11.00) {};
				\vertex (j) at (4.1,-11.00) {};
				\vertex (k) at (0,24.00) {};
				\draw[edge] (a)edge(e) (a)edge(f) (a)edge(g) (a)edge(h) (a)edge(k)
					(b)edge(c) (b)edge(e) (b)edge(g) (b)edge(i) (c)edge(f) (c)edge(h)
					(c)edge(j) (d)edge(e) (d)edge(f) (d)edge(i) (d)edge(j) (g)edge(k)
					(h)edge(k) (i)edge(j);
			\end{tikzpicture}
			&
			\begin{tikzpicture}[scale=0.06]
				\vertex (1) at (10.4973, 0.567839) {};
				\vertex (2) at (13.099, -10.2247) {};
				\vertex (3) at (4.22603, 6.35668) {};
				\vertex (4) at (-10.1444, -3.38527) {};
				\vertex (5) at (-0.546292, 18.3514) {};
				\vertex (6) at (-7.75724, -13.) {};
				\vertex (7) at (20.3634, 0.567839) {};
				\vertex (8) at (4.22603, -5.221) {};
				\vertex (9) at (-10.1444, 4.52095) {};
				\vertex (10) at (0.49548, 0.330676) {};
				\vertex (11) at (-4.67816, 0.567839) {};
				\vertex (12) at (-19.6366, 0.567839) {};

				\draw[edge] (1)edge(2) (1)edge(3) (1)edge(7) (1)edge(8) (2)edge(4)
				(2)edge(7) (2)edge(8) (3)edge(5) (3)edge(9) (3)edge(10) (4)edge(6)
				(4)edge(9) (4)edge(11) (5)edge(7) (5)edge(11) (5)edge(12) (6)edge(8)
				(6)edge(10) (6)edge(12) (9)edge(12) (10)edge(11) ;
			\end{tikzpicture}
			\\
			688 & 1760 & 4864\\
		\end{tabular}
	\end{center}
	\caption{Laman graphs with $7\leq n\leq12$ vertices that have the 4-vertex
		Laman graph (encoded as 31) as a subgraph; below their Laman numbers are given. In some
		cases there are several Laman graphs with this subgraph property and with
		the same Laman number, but among all Laman graphs that have this subgraph
		there does not exist one with higher Laman number (encodings see Table~\ref{table:enc:31-fan}).}
	\label{figure:sub4_max}
\end{figure}
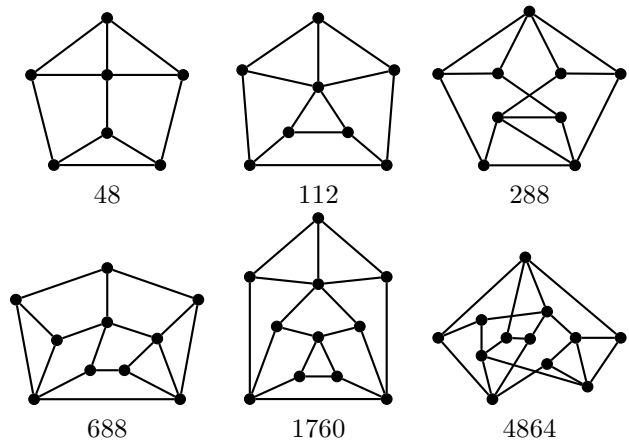

\begin{table}[H]
	\begin{center}
		\begin{tabular}{rllllll}
			\toprule
			$n$ & caterpillar & fan $T(n)$ & fan     & 31-fan  & 254-fan & 7916-fan\\
			\midrule
			6   & 2.21336     & 2.28943    &         & 2       & 2       & -       \\
			7   & 2.23685     & 2.30033    &         & 2.28943 & 2       & 2       \\
			8   & 2.26772     & 2.32542    &         & 2.30033 & 2.28943 & 2       \\
			9   & 2.30338     & 2.35824    &         & 2.35216 & 2.30033 & 2.28943 \\
			10  & 2.33378     & 2.38581    &         & 2.35824 & 2.35216 & 2.30033 \\
			11  & 2.36196     & 2.41159    &         & 2.38581 & 2.35824 & 2.35216 \\
			12  & 2.39386     & 2.43198    &         & 2.43006 & 2.39802 & 2.35824 \\\hline
			13  & 2.40453     & 2.44156    & 2.44498 & 2.44772 & 2.42197 & 2.39802 \\
			14  & 2.43185     & 2.45868    & 2.46087 & 2.46391 & 2.44251 & 2.42197 \\
			15  & 2.44695     & 2.47445    &         & 2.47076 & 2.45031 & 2.42906 \\
			16  & 2.46890     & 2.48657    &         & 2.48794 & 2.47166 & 2.43712 \\
			17  & 2.48875     & 2.49668    & 2.49779 & 2.49160 & 2.48043 & 2.46341 \\
			18  & 2.49378     & 2.50798    & \\
			\bottomrule
		\end{tabular}
	\end{center}
	\caption{Growth rates (rounded) of the lower bounds.
					For $n\leq12$ these values
          are proven to be the best achievable ones; for $n>12$ the values are
          just the best we found by experiments, hence it is possible that
          there are better ones.
          The drawings of the graphs corresponding to
          the last three columns are given in Figure~\ref{figure:fanbases}.\\
          The encodings for the graphs can be found at: caterpillar (Table~\ref{table:enc:max_laman_numbers}),
          fan $T(n)$ (Table~\ref{table:enc:max_graphs_triangle}), fan (Table~\ref{table:enc:fan}),
          31-fan (Table~\ref{table:enc:31-fan}), 254-fan (Table~\ref{table:enc:254-fan}), 7916-fan (Table~\ref{table:enc:7916-fan})
          }
	\label{table:bounds}
\end{table}

\begin{figure}[H]
\begin{center}
\begin{tikzpicture}[scale=1.15,value/.style={fill=black,circle,inner sep=1pt}]
	\begin{scope}[yscale=20,xscale=0.7]
    \draw[->] (5,2.19) -- (5,2.55) node[left,align=left] {growth\\rate};
    \draw[->] (5,2.19) -- (19,2.19) node[below] {$n$};
    \foreach \y in {2.2,2.3,2.4,2.5} \draw (5,\y) -- +(-0.1,0) node[left] {$\y$};
    \foreach \y in {2.2,2.21,...,2.52} \draw (5,\y) -- +(-0.05,0);
    \foreach \x in {6,...,18} \draw (\x,2.19) -- +(0,-0.005) node[below] {$\x$};
    \draw[dashed,black!40!white] (12,2.19) -- (12,2.55);
		\draw[Red, line width=1pt]
			 (6, 2.21336)
		-- (7, 2.23685)
		-- (8, 2.26772)
		-- (9, 2.30338)
		-- (10, 2.33378)
		-- (11, 2.36196)
		-- (12, 2.39386) ;
		\draw[Red!30!white,line width=1pt]
		   (12, 2.39386)
		-- (13, 2.40453)
		-- (14, 2.43185)
		-- (15, 2.44695)
		-- (16, 2.4689)
		-- (17, 2.48875)
		-- (18,2.49378);
		\node[value] at (6, 2.21336) {};
		\node[value] at (7, 2.23685)  {};
		\node[value] at (8, 2.26772)  {};
		\node[value] at (9, 2.30338)  {};
		\node[value] at (10, 2.33378)  {};
		\node[value] at (11, 2.36196)  {};
		\node[value] at (12, 2.39386)  {};
		\node[value] at (13, 2.40453)  {};
		\node[value] at (14, 2.43185)  {};
		\node[value] at (15, 2.44695)  {};
		\node[value] at (16, 2.4689)  {};
		\node[value] at (17, 2.48875)  {};
		\node[value] at (18, 2.49378) {};
		\draw[Blue, line width=1pt]
			 (6, 2.28943)
		-- (7, 2.30033)
		-- (8, 2.32542)
		-- (9, 2.35824)
		-- (10, 2.38581)
		-- (11, 2.41159)
		-- (12, 2.43198);
		\draw[Blue!30!white, line width=1pt]
			 (12, 2.43198)
		-- (13, 2.44498)
		-- (14, 2.46087)
		-- (15, 2.47445)
		-- (16, 2.48657)
		-- (17, 2.49779)
		-- (18, 2.50798);
		\node[value] at (6, 2.28943) {};
		\node[value] at (7, 2.30033)  {};
		\node[value] at (8, 2.32542)  {};
		\node[value] at (9, 2.35824)  {};
		\node[value] at (10, 2.38581)  {};
		\node[value] at (11, 2.41159)  {};
		\node[value] at (12, 2.43198)  {};
		\node[value] at (13, 2.44498)  {};
		\node[value] at (14, 2.46087)  {};
		\node[value] at (15, 2.47445)  {};
		\node[value] at (16, 2.48657)  {};
		\node[value] at (17, 2.49779)  {};
		\node[value] at (18, 2.50798)  {};
		\draw[Green, line width=1pt]
			 (7, 2.28943)
		-- (8, 2.30033)
		-- (9, 2.35216)
		-- (10, 2.35824)
		-- (11, 2.38581)
		-- (12, 2.43006);
		\draw[Green!30!white, line width=1pt]
		   (12,2.43006)
		-- (13,2.44772)
		-- (14,2.46391)
		-- (15,2.47076)
		-- (16,2.48794)
		-- (17,2.49160);
		\node[value] at (7, 2.28943)  {};
		\node[value] at (8, 2.30033)  {};
		\node[value] at (9, 2.35216)  {};
		\node[value] at (10, 2.35824)  {};
		\node[value] at (11, 2.38581)  {};
		\node[value] at (12, 2.43006) {};
		\node[value] at (13,2.44772) {};
		\node[value] at (14,2.46391) {};
		\node[value] at (15,2.47076) {};
		\node[value] at (16,2.48794) {};
		\node[value] at (17,2.49160) {};
		\draw[Fuchsia, line width=1pt]
			 (8, 2.28943)
		-- (9, 2.30033)
		-- (10, 2.35216)
		-- (11, 2.35824)
		-- (12, 2.39802);
		\draw[Fuchsia!30!white, line width=1pt]
		   (12,2.39802)
		-- (13,2.42197)
		-- (14,2.44251)
		-- (15,2.45031)
		-- (16,2.47166)
		-- (17,2.48043);
		\node[value] at (8, 2.28943)  {};
		\node[value] at (9, 2.30033)  {};
		\node[value] at (10, 2.35216) {};
		\node[value] at (11, 2.35824) {};
		\node[value] at (12, 2.39802) {};
		\node[value] at (13, 2.42197) {};
		\node[value] at (14, 2.44251) {};
		\node[value] at (15, 2.45031) {};
		\node[value] at (16, 2.47166) {};
		\node[value] at (17, 2.48043) {};
		\draw[Brown, line width=1pt]
			 (9, 2.28943)
		-- (10, 2.30033)
		-- (11, 2.35216)
		-- (12, 2.35824);
		\draw[Brown!30!white, line width=1pt]
		   (12,2.35824)
		-- (13,2.39802)
		-- (14,2.42197)
		-- (15,2.42906)
		-- (16,2.43712)
		-- (17,2.46341);
		\node[value] at (9, 2.28943)  {};
		\node[value] at (10, 2.30033)  {};
		\node[value] at (11, 2.35216)  {};
		\node[value] at (12, 2.35824) {};
		\node[value] at (13,2.39802) {};
		\node[value] at (14,2.42197) {};
		\node[value] at (15,2.42906) {};
		\node[value] at (16,2.43712) {};
		\node[value] at (17,2.46341) {};
  \end{scope}
\end{tikzpicture}
\end{center}
\caption{Growth rates of the lower bounds (red = caterpillar construction,
  blue = fan construction, green = 31-fan construction, fuchsia = 254-fan construction, brown = 7916-fan construction).
  The light colors indicate values that were not found by exhaustive search and which therefore could possibly
  be improved.}
\label{fig:growth}
\end{figure}
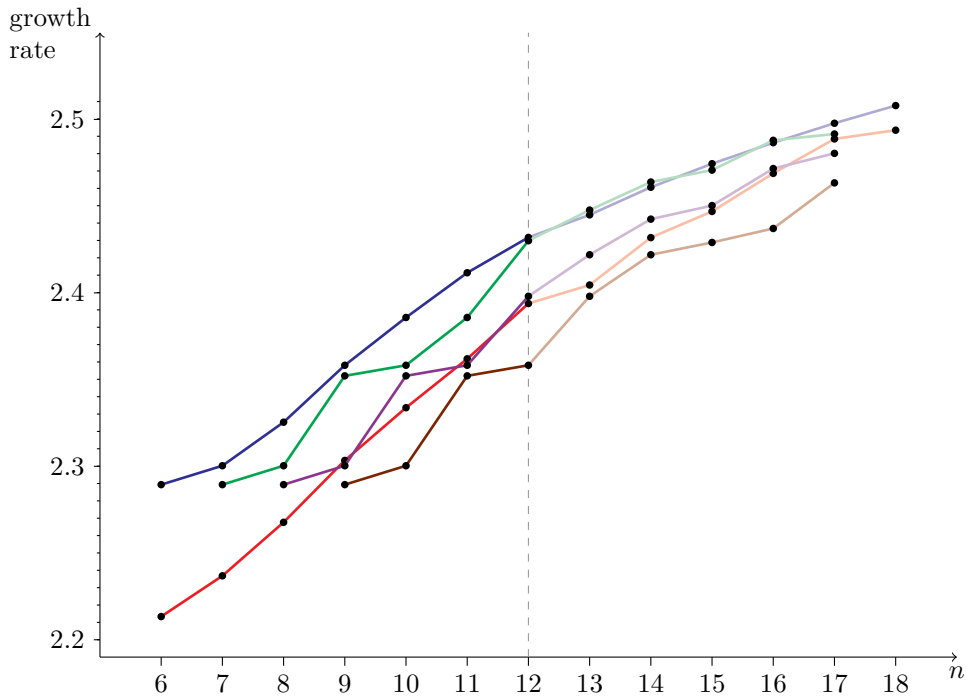

From Table~\ref{table:bounds} we can see the bound $2.28943^n$
obtained in~\cite{Borcea2004}, $2.30033^n$ from~\cite{EmirisMoroz},
and $2.41159^n$ from~\cite{Jackson2018}, as well as the current
improvements obtained in this paper.
By instantiating Formula~\eqref{eq:bound_fan} with the last Laman graph in Figure~\ref{figure:max_graphs_triangle}, which has 18 vertices and Laman number $1953816$,
we obtain the following theorem.
\begin{theorem}
  \label{thm:lower-bd-2d}
  The maximal Laman number $\maxlamIIn$ satisfies
  \begin{equation*}
		\maxlamIIn \geq  2\cdot2^{(n-3)\modop 15}\cdot976908^{\lfloor(n-3)/15\rfloor}.
	\end{equation*}
  This means $\maxlamIIn$ grows at least as $\bigl(\!\sqrt[15]{976908}\bigr)^n$, which is approximately $2.50798^n$.
  In other words $\bigl(\!\sqrt[15]{976908}\bigr)^n\in\mathcal O(\maxlamIIn)$.
\end{theorem}

\section{Dimension 3}
\label{sec:dim-3}
A generalization of the counting condition to three dimensions would suggest that
a graph $G=(V,E)$ needs to fulfill $\vert E\vert = 3\vert
V\vert-6$, and $\vert E' \vert \leq 3\vert V' \vert-6$ for
every subgraph $G'=(V',E')$ of $G$ in order to be rigid.
Unlike the two-dimensional case, this definition is necessary but not
sufficient for generic minimal rigidity.  An example of a graph
which is not minimally rigid in dimension~$3$ can already be found
in~\cite{Geiringer1927}.  We are interested in lower bounds on
$\maxlamIIIn$, which is the three-dimensional analog of $\maxlamIIn$.

\begin{definition}
  Let $G=(V,E)$ be a graph.  We call $G$ a Geiringer graph%
  \footnote{As Hilda Pollaczek-Geiringer had already worked on rigid
      graphs in 2D and 3D \cite{Geiringer1927,Geiringer1932}, long before
      Gerard Laman~\cite{Laman1970}.}, if there exists only a finite number of
    (complex) spatial embeddings in $\C^3$, given a generic labeling
    $\lambda\colon E\to\C$ of the edges of~$G$.

  For a Geiringer graph $G$, we define $\lamIII(G)$, called
  the \emph{3D-Laman number} of $G$, to be this finite number of
  (complex) embeddings. Moreover, we define $\maxlamIIIn$ to be the
  largest 3D-Laman number that is achieved among all Geiringer graphs with $n$ vertices.
\end{definition}

\begin{figure}[H]
  \begin{center}
		\begin{subfigure}[b]{0.98\textwidth}
      \begin{center}
				\begin{tikzpicture}[scale=1.7]
					\vertex (a) at (0,0) {};
					\vertex (b) at (1,0) {};
					\vertex (c) at (0.5,0.866025) {};
					\vertex (d) at (2.5,0) {};
					\vertex (e) at (3.5,0) {};
					\vertex (f) at (3,0.866025) {};
					\node[nvertex] (g) at (3,0.33) {};
					\draw[edgeq] (a)edge(b) (a)edge(c) (b)edge(c) (d)edge(e) (d)edge(f) (e)edge(f);
					\draw[nedge] (d)edge(g) (e)edge(g) (f)edge(g);
					\draw[ultra thick,->] (1.5,0.5) -- (2,0.5);
				\end{tikzpicture}
				\caption{Type~1}\label{fig:3dtype1}
      \end{center}
    \end{subfigure}
		\\
		\begin{subfigure}[b]{0.98\textwidth}
      \begin{center}
				\begin{tikzpicture}[scale=1.7]
					\vertex (a) at (0,0) {};
					\vertex (b) at (1,0) {};
					\vertex (c) at (1,1) {};
					\vertex (d) at (0,1) {};
					\vertex (e) at (2.5,0) {};
					\vertex (f) at (3.5,0) {};
					\vertex (g) at (3.5,1) {};
					\vertex (h) at (2.5,1) {};
					\node[nvertex] (i) at (3,0.25) {};
					\draw[oedge] (a)edge(b);
					\draw[edgeq] (a)edge(c) (a)edge(d) (b)edge(c) (b)edge(d) (c)edge(d);
					\draw[edgeq] (e)edge(g) (e)edge(h) (f)edge(g) (f)edge(h) (g)edge(h);
					\draw[nedge]  (e)edge(i) (f)edge(i) (g)edge(i) (h)edge(i);
					\draw[ultra thick,->] (1.5,0.5) -- (2,0.5);
				\end{tikzpicture}
				\caption{Type~2}\label{fig:3dtype2}
      \end{center}
    \end{subfigure}
    \\
    \begin{subfigure}[b]{0.98\textwidth}
      \begin{center}
				\begin{tikzpicture}[scale=1.7]
					\begin{scope}[scale=0.75]
						\vertex (a) at (0.587785, -0.809017) {};
						\vertex (b) at (0.951057, 0.309017) {};
						\vertex (c) at (0., 1.) {};
						\vertex (d) at (-0.951057,0.309017) {};
						\vertex (e) at (-0.587785, -0.809017) {};
						
						\draw[edgeq] (a)edge(b) (a)edge(c) (a)edge(e) (b)edge(c) (b)edge(d) (c)edge(d) (c)edge(e) (d)edge(e);
						\draw[oedge] (a)edge(d) (b)edge(e);
					\end{scope}
					\draw[ultra thick,->] (1.5,0.) -- (2,0.);
					\begin{scope}[xshift=3.5cm,scale=0.75]
						\vertex (a) at (0.587785, -0.809017) {};
						\vertex (b) at (0.951057, 0.309017) {};
						\vertex (c) at (0., 1.) {};
						\vertex (d) at (-0.951057,0.309017) {};
						\vertex (e) at (-0.587785, -0.809017) {};
						\node[nvertex] (f) at (0,0) {};
						\draw[edgeq] (a)edge(b) (a)edge(c) (a)edge(e) (b)edge(c) (b)edge(d) (c)edge(d) (c)edge(e) (d)edge(e);
						\draw[nedge] (a)edge(f) (b)edge(f) (c)edge(f) (d)edge(f) (e)edge(f);
					\end{scope}
				\end{tikzpicture}
				\caption{Type~3x}\label{fig:3dtype3x}
      \end{center}
    \end{subfigure}
    \\
    \begin{subfigure}[b]{0.98\textwidth}
      \begin{center}
				\begin{tikzpicture}[scale=1.7]
					\begin{scope}[scale=0.75]
						\vertex (a) at (0.587785, -0.809017) {};
						\vertex (b) at (0.951057, 0.309017) {};
						\vertex (c) at (0., 1.) {};
						\vertex (d) at (-0.951057,0.309017) {};
						\vertex (e) at (-0.587785, -0.809017) {};
						
						\draw[edgeq] (a)edge(b) (a)edge(d) (a)edge(e) (b)edge(c) (b)edge(d) (c)edge(d) (b)edge(e) (d)edge(e);
						\draw[oedge] (a)edge(c) (c)edge(e);
					\end{scope}
					\draw[ultra thick,->] (1.5,0.) -- (2,0.);
					\begin{scope}[xshift=3.5cm,scale=0.75]
						\vertex (a) at (0.587785, -0.809017) {};
						\vertex (b) at (0.951057, 0.309017) {};
						\vertex (c) at (0., 1.) {};
						\vertex (d) at (-0.951057,0.309017) {};
						\vertex (e) at (-0.587785, -0.809017) {};
						\node[nvertex] (f) at (0,0) {};
						\draw[edgeq] (a)edge(b) (a)edge(d) (a)edge(e) (b)edge(c) (b)edge(d) (c)edge(d) (b)edge(e) (d)edge(e);
						\draw[nedge] (a)edge(f) (b)edge(f) (c)edge(f) (d)edge(f) (e)edge(f);
					\end{scope}
				\end{tikzpicture}
				\caption{Type~3v}\label{fig:3dtype3v}
      \end{center}
    \end{subfigure}
  \end{center}
  \caption{Henneberg steps of different types in dimension 3; a dashed line indicates that this edge can exist but does not need to.}
  \label{fig:Henneberg3d}
\end{figure}
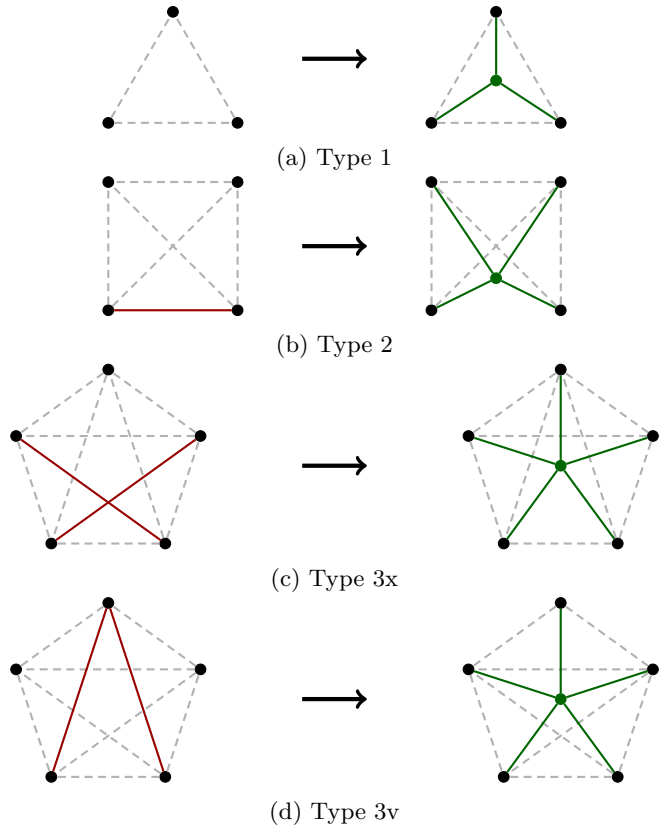

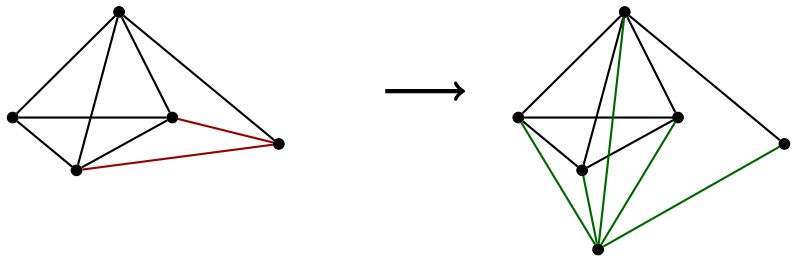
\begin{figure}[H]
	\begin{center}
		\begin{tikzpicture}[scale=0.7]
			\vertex (a) at (2,3) {};
			\vertex (b) at (1.2,0) {};
			\vertex (c) at (0,1) {};
			\vertex (d) at (3,1) {};
			\vertex (e) at (5,0.5) {};
			\draw[edge] (a)edge(b) (a)edge(c) (a)edge(d) (a)edge(e) (b)edge(c) (b)edge(d) (c)edge(d);
			\draw[oedge] (b)edge(e) (d)edge(e);
			\draw[ultra thick,->] (7,1.5) -- (8.5,1.5);
			\begin{scope}[xshift=9.5cm]
				\vertex (a) at (2,3) {};
				\vertex (b) at (1.2,0) {};
				\vertex (c) at (0,1) {};
				\vertex (d) at (3,1) {};
				\vertex (e) at (5,0.5) {};
				\vertex (f) at (1.5,-1.5) {};
				\draw[edge] (a)edge(b) (a)edge(c) (a)edge(d) (a)edge(e) (b)edge(c) (b)edge(d) (c)edge(d);
				\draw[nedge] (a)edge(f) (b)edge(f) (c)edge(f) (d)edge(f) (e)edge(f);
			\end{scope}
		\end{tikzpicture}
	\end{center}
	\caption{Flexible graph constructed by a Henneberg move of type 3}
	\label{fig:flex3v}
\end{figure}

In \cite{TayWhiteley} Geiringer graphs are shown to be constructible from a triangle graph
by a sequence of three types of steps (see Figure~\ref{fig:Henneberg3d}).
Steps of type~1 and type~2 preserve rigidity (see \cite{TayWhiteley}).
The steps of type 3 can be further classified according to whether the two chosen edges have a common vertex or not.
Note that every Geiringer graph can be constructed using such steps \cite[Prop.~4.1, 4.4, 4.5]{TayWhiteley}, but not every construction by these steps is indeed minimally rigid,
i.e.\ rigidity is not necessarily preserved by steps of type 3. Indeed type 3v does not even preserve the vertex-edge-count for subgraphs (see Figure~\ref{fig:flex3v}).
However, there are certain subclasses of type 3 steps for which rigidity is preserved (see for instance \cite{TayWhiteley,GraverServatius,Cruickshank2014}).

In the following we construct Geiringer graphs by the above mentioned moves, removing those which turn out to be non-rigid.
By this procedure we get all Geiringer graphs with up to 10 vertices.
The computation of the number of realizations is done by Gr\"obner bases:
The coordinates of the vertices are obtained as the solutions of a system
of (quadratic) polynomial equations. Instead of keeping the edge lengths generic
(by introducing a symbolic parameter for each edge), we insert random numbers (integers) for
the edge lengths. Otherwise the computation would not be feasible at all.
Moreover, for further speed-up, we compute the Gr\"obner basis only modulo a sufficiently large prime
number~$p$ so that the occurrence of large rational numbers is avoided.
In other words, the Gr\"obner basis computation takes place over the finite field~$\Z_p$.
In order to get high confidence into the results, we did each computation at least
three times, with different random choices of the parameters. If we get the same result
three times, we can be rather sure to have the correct number. However, we want to
make the reader aware of the fact, that it is a probabilistic method.
Although we have a strong evidence for the computed 3D-Laman numbers, they are not
rigorously proven to be correct.

Still, computing the 3D-Laman numbers for all Geiringer graphs of 10
vertices was a major undertaking. By applying the Henneberg steps depicted in
Figure~\ref{fig:Henneberg3d} in all possible ways, we obtained 747065 graphs
that potentially had the property of being minimally rigid (our Gr\"obner
basis computations suggested that 612884 of them indeed have this
property). In our implementation, we do some preprocessing on the graphs in
order to create polynomial systems with as few variables as possible: for
example, we remove vertices of valency~$3$ (i.e., revert Henneberg steps of
type~$1$), and compensate by multiplying the final Laman number by~$2$ for
each removed vertex. Another optimization consists in identifying the largest
tetrahedral subgraph, i.e., the largest subgraph that can be constructed by
Henneberg steps of type~$1$, starting from a triangle. This subgraph is
considered when fixing some vertices of the graph, in order to deal with
rotations and translations.  Then we call the fast FGb~\cite{FGb}
implementation of Gr\"obner bases in Maple, for determining the number of
solutions of the constructed polynomial system. Executing this program once
for all 747065 graphs took about 162 days of CPU time, using Xeon E5-2630v3
Haswell 2,4Ghz CPUs.  However, the computations were run in parallel so that
the result was obtained after a few days. This means that in average it took
about 19s to determine the 3D-Laman number of a graph with 10 vertices, but
the timings vary a lot: graphs which can be constructed by Henneberg steps of
type~1 require almost no time, due to our preprocessing, while the Gr\"obner
basis computation for some graphs takes several hours (up to 16 hours).

It is easy to see that a Henneberg step of type~1 always increases the 3D-Laman number by a factor of~$2$.
So far it is not known by which factor a Henneberg step of type~2 or type~3 might increase the 3D-Laman number.
Table~\ref{table:Laman_Increase_3d} summarizes some increases of 3D-Laman numbers, given a certain Geiringer graph $G$ and constructing a new one $G'$ by a single Henneberg step.
\begin{table}[H]
  \begin{center}
    \begin{tabular}{lrrrrr}
      \toprule
      Type & $G$         & $\lamIII(G)$ & $G'$          & $\lamIII(G')$ & Factor\\\midrule
      3v   & 11717490611 & 512    & 9634462543324 & 128     & 0.25\\
      3v   & 49724126    & 160    & 18848282483   & 64      & 0.40\\
      3v   & 515806      & 48     & 203906043     & 32      & 0.66\\
      2    & 981215      & 24     & 31965132      & 24      & 1.00\\
      3x   & 16350       & 16     & 1973983       & 16      & 1.00\\      
      2, 3x & 1973983     & 16     & 49524604      & 128     & 8.00\\
      3x   & 384510      & 16     & 49724126      & 160     & 10.00\\
      3v   & 382463      & 16     & 49724126      & 160     & 10.00\\
      3x   & 15661790    & 32     & 7309884067    & 512     & 16.00\\
      3x   & 2000476603  & 48     & 2704137746603 & 1088    & 22.66\\
      \bottomrule
    \end{tabular}
  \end{center}
  \caption{Henneberg constructions and increase of 3D-Laman numbers}
  \label{table:Laman_Increase_3d}
\end{table}

\subsection{Constructions}
\label{sec:constr-3d}

We consider again caterpillar and fan constructions.
For the caterpillar we now need to glue two graphs by a common triangle.
Similarly, we need a tetrahedron for the fan construction.
For the generalized fan construction we use the unique Geiringer graph with $5$ vertices.
For sake of completeness, we display the general formula for obtaining a lower bound on
$\maxlamIIIn$ with the generalized 3D-fan construction; the formula is completely analogous
to~\eqref{eq:bound_genfan}:
\begin{equation}\label{eq:bound_genfan_3d}
  \maxlamIIIn \geq 2^{(n-|W|)\modop(|V|-|W|)} \cdot\lamIII(H) \cdot
  \left(\frac{\lamIII(G)}{\lamIII(H)}\right)^{\!\lfloor(n-|W|)/(|V|-|W|)\rfloor} \qquad (n\geq|W|).
\end{equation}

\subsection[Lower Bounds for M3(n)]{Lower Bounds for $\maxlamIIIn$}
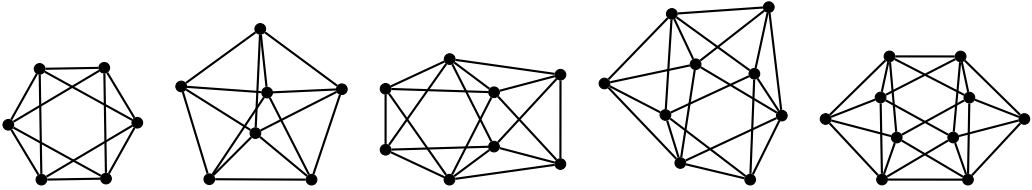
\begin{figure}[H]
  \begin{center}
    \begin{tabular}{c@{\hspace{12pt}}c@{\hspace{12pt}}c@{\hspace{12pt}}c@{\hspace{12pt}}c}
      \begin{tikzpicture}[scale=1.2]
        \node[vertex] (1) at (1.07046, 0.0117707) {};
				\node[vertex] (2) at (0.342906, 1.22397) {};
				\node[vertex] (3) at (1.41347, 0.629078) {};
				\node[vertex] (4) at (0., 0.606446) {};
				\node[vertex] (5) at (0.362485, 0.) {};
				\node[vertex] (6) at (1.05135, 1.23624) {};
				\draw[edge] (1)edge(3) (1)edge(4) (1)edge(5) (1)edge(6) (2)edge(3) 
				(2)edge(4) (2)edge(5) (2)edge(6) (3)edge(5) (3)edge(6) (4)edge(5) 
				(4)edge(6) ;
      \end{tikzpicture}
      &
      \begin{tikzpicture}[scale=1.2]
        \node[vertex] (1) at (1.42865, 0.) {};
				\node[vertex] (2) at (0.866612, 1.66585) {};
				\node[vertex] (3) at (0., 1.02986) {};
				\node[vertex] (4) at (0.308272, 0.00486401) {};
				\node[vertex] (5) at (1.76067, 0.999258) {};
				\node[vertex] (6) at (0.941843, 0.96197) {};
				\node[vertex] (7) at (0.813498, 0.514202) {};
				\draw[edge] (1)edge(4) (1)edge(5) (1)edge(6) (1)edge(7) (2)edge(3) 
				(2)edge(5) (2)edge(6) (2)edge(7) (3)edge(4) (3)edge(6) (3)edge(7) 
				(4)edge(6) (4)edge(7) (5)edge(6) (5)edge(7) ;
      \end{tikzpicture}
      &
      \begin{tikzpicture}[scale=1.2]
        \node[vertex] (1) at (1.91578, 0.172323) {};
				\node[vertex] (2) at (0.00021641, 1.00356) {};
				\node[vertex] (3) at (0., 0.33047) {};
				\node[vertex] (4) at (1.9166, 1.15968) {};
				\node[vertex] (5) at (0.702382, 1.33239) {};
				\node[vertex] (6) at (0.699981, 0.) {};
				\node[vertex] (7) at (1.18809, 0.365389) {};
				\node[vertex] (8) at (1.18866, 0.965653) {};
				\draw[edge] (1)edge(4) (1)edge(6) (1)edge(7) (1)edge(8) (2)edge(3) 
				(2)edge(5) (2)edge(6) (2)edge(8) (3)edge(5) (3)edge(6) (3)edge(7) 
				(4)edge(5) (4)edge(7) (4)edge(8) (5)edge(7) (5)edge(8) (6)edge(7) 
				(6)edge(8) ;
      \end{tikzpicture}
      &
      \begin{tikzpicture}[scale=1.2]
        \node[vertex] (1) at (0., 1.06259) {};
				\node[vertex] (2) at (1.59714, 0.) {};
				\node[vertex] (3) at (1.801, 1.90641) {};
				\node[vertex] (4) at (1.64311, 1.17024) {};
				\node[vertex] (5) at (0.998848, 1.27554) {};
				\node[vertex] (6) at (0.83117, 0.184029) {};
				\node[vertex] (7) at (1.94377, 0.707595) {};
				\node[vertex] (8) at (0.668086, 0.714625) {};
				\node[vertex] (9) at (0.738407, 1.83153) {};
				\draw[edge] (1)edge(5) (1)edge(6) (1)edge(8) (1)edge(9) (2)edge(4) 
				(2)edge(6) (2)edge(7) (2)edge(8) (3)edge(4) (3)edge(5) (3)edge(7) 
				(3)edge(9) (4)edge(7) (4)edge(8) (4)edge(9) (5)edge(6) (5)edge(7) 
				(5)edge(9) (6)edge(7) (6)edge(8) (8)edge(9) ;
      \end{tikzpicture}
      &
      \begin{tikzpicture}
        \node[vertex] (1) at (2.61405, 0.804397) {};
				\node[vertex] (2) at (0., 0.804295) {};
				\node[vertex] (3) at (0.840417, 1.63439) {};
				\node[vertex] (4) at (0.742173, 0.00138984) {};
				\node[vertex] (5) at (1.87222, 0.) {};
				\node[vertex] (6) at (1.67869, 0.559079) {};
				\node[vertex] (7) at (0.723064, 1.08966) {};
				\node[vertex] (8) at (0.935239, 0.558165) {};
				\node[vertex] (9) at (1.77585, 1.63352) {};
				\node[vertex] (10) at (1.89103, 1.0882) {};
				\draw[edge] (1)edge(5) (1)edge(6) (1)edge(9) (1)edge(10) (2)edge(3)
				(2)edge(4) (2)edge(7) (2)edge(8) (3)edge(7) (3)edge(8) (3)edge(9)
				(3)edge(10) (4)edge(5) (4)edge(6) (4)edge(7) (4)edge(8) (5)edge(6)
				(5)edge(8) (5)edge(10) (6)edge(7) (6)edge(9) (7)edge(9) (8)edge(10)
				(9)edge(10) ;
      \end{tikzpicture}
    \end{tabular}
  \end{center}
  \caption{Geiringer graphs with $6\leq n\leq10$ vertices; for each $n$ the (unique)
  graph with maximal number of embeddings is depicted. The corresponding 3D-Laman
  numbers $\maxlamIIIn$ are given in Table~\ref{table:max_laman_numbers_3d} (encodings see Table~\ref{table:enc:max_graphs_3d}).}
  \label{figure:max_graphs_3d}
\end{figure}

In order to get good lower bounds, we need particular Geiringer graphs that have a
large number of embeddings. We computed the 3D-Laman numbers of all Geiringer graphs
with up to $n=10$ vertices. For each $n$ we have identified the (unique) Geiringer graph
with the highest number of embeddings. These numbers are given in
Table~\ref{table:max_laman_numbers_3d}. The corresponding graphs for $6\leq
n\leq10$ are shown in Figure~\ref{figure:max_graphs_3d}.
\begin{table}[H]
	\begin{center}
		\begin{tabular}{lcccccccc}
			\toprule
			$n$           & 6  & 7   & 8    & 9    & 10 & 11 & 12 \\
			min           & 8  & 16  & 24   & 48   & 76 \\
			$\maxlamIIIn$ & 16 & 48  & 160  & 640  & 2560\\
			upper         & 40 & 224 & 1344 & 8448 & 54912 & 366080 & 2489344 \\
			\bottomrule
		\end{tabular}
	\end{center}
	\caption{Minimal and maximal 3D-Laman number among all $n$\hbox{-}vertex Geiringer graphs;
	the row labeled with ``min'' contains the lowest 3D-Laman number which is found by computation.
	The row labeled with ``upper'' contains the bounds from \cite{Borcea2004}.}
	\label{table:max_laman_numbers_3d}
\end{table}

In \cite{Emiris2009} lower and upper bounds for the 1-skeleta of simplicial polyhedra are computed.
They also use an extension of Henneberg steps to the three-dimensional case. However,
they form just a subset of the Henneberg steps presented here.
From Table~\ref{table:bounds3d} we can see the bound of $2.51984^n$
obtained in \cite{Emiris2009} and the improvements obtained in this
paper.
By instantiating Formula~\eqref{eq:bound_genfan_3d} with the last Laman graph in Figure~\ref{figure:max_graphs_3d}, which has 10 vertices and 3D-Laman number $2560$,
we obtain the following theorem.
\begin{theorem}
  \label{thm:lower-bd-3d}
  The maximal Laman number $\maxlamIIIn$ satisfies
  \begin{equation*}
		\maxlamIIIn \geq 2^{(n-3)\modop 7} \cdot 2560^{\lfloor(n-3)/7\rfloor}.
	\end{equation*}
  This means $\maxlamIIIn$ grows at least as $\bigl(\!\sqrt[7]{2560}\bigr)^n$ which is approximately $3.06825^n$.
  In other words $\bigl(\!\sqrt[7]{2560}\bigr)^n\in\mathcal O(\maxlamIIIn)$.
\end{theorem}

\begin{table}[H]
	\begin{center}
		\begin{tabular}{rlll}
			\toprule
			$n$ & caterpillar & fan     & generalized fan \\\midrule
			6   & 2.51984     & 2       & -\\
			7   & 2.63215     & 2.51984 & 2 \\
			8   & 2.75946     & 2.63215 & 2.51984\\
			9   & 2.93560      & 2.95155 & 2.82843\\
			10  & 3.06825     & 3.06681 & 2.95155\\
			\bottomrule
		\end{tabular}
	\end{center}
	\caption{Growth rates (rounded) of the lower bounds. The encodings for the graphs can be found at:
	caterpillar (Table~\ref{table:enc:max_graphs_3d}), fan (Table~\ref{table:enc:fan_3d}), generalized fan (Table~\ref{table:enc:genfan_3d})}
	\label{table:bounds3d}
\end{table}

\begin{figure}[H]
\begin{center}
\begin{tikzpicture}[scale=1,value/.style={fill=black,circle,inner sep=1pt}]
  \begin{scope}[xscale=1.3,yscale=8]
    \draw[->] (5,2.4) -- (5,3.2) node[left] {growth rate};
    \draw[->] (5,2.4) -- (11,2.4) node[below] {$n$};
    \foreach \y in {2.5,2.6,2.7,2.8,2.9,3.0,3.1} \draw (5,\y) -- +(-0.1,0) node[left] {$\y$};
    \foreach \x in {6,...,10} \draw (\x,2.4) -- +(0,-0.03) node[below] {$\x$};
		\draw[Red, line width=1pt]
			 (6, 2.51984)
		-- (7, 2.63215)
		-- (8, 2.75946)
		-- (9, 2.9356)
		-- (10, 3.06825);
		\node[value] at (6, 2.51984) {};
		\node[value] at (7,2.63215)  {};
		\node[value] at (8, 2.75946)  {};
		\node[value] at (9, 2.9356)  {};
		\node[value] at (10, 3.06825) {};
		\draw[Blue, line width=1pt]
			 (7, 2.51984) -- (8, 2.63215) -- (9, 2.95155) -- (10,3.06681);
		\node[value] at (7, 2.51984) {};
		\node[value] at (8, 2.63215) {};
		\node[value] at (9, 2.95155) {};
		\node[value] at (10,3.06681) {};
		\draw[Green, line width=1pt]
			 (8, 2.51984) -- (9, 2.82843) -- (10,2.95155);
		\node[value] at (8, 2.51984) {};
		\node[value] at (9, 2.82843) {};
		\node[value] at (10,2.95155) {};
  \end{scope}
\end{tikzpicture}
\end{center}
\caption{Growth rates of the lower bounds (red = caterpillar construction,
  blue = fan construction, green = generalized fan construction)}
\label{fig:growth_3d}
\end{figure}
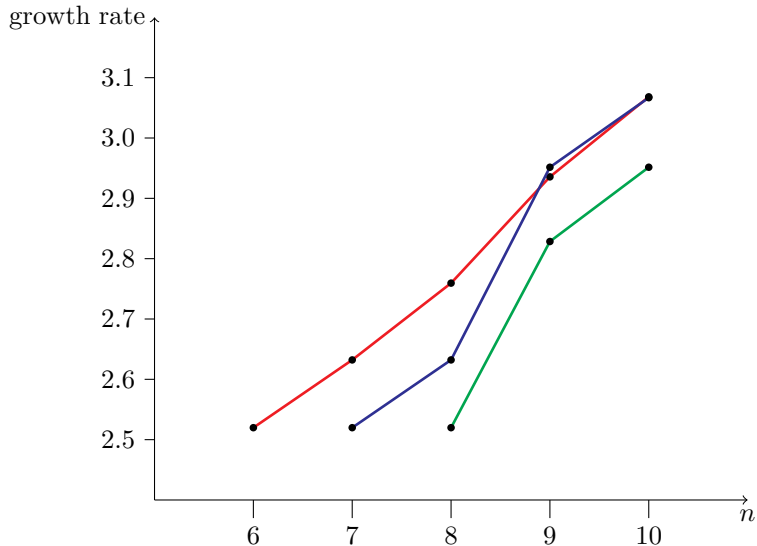

\section{Conclusion}
By exploiting  state-of-the-art methods we gave some new bounds on the maximal
possible number of realizations of rigid graphs for a given number of vertices.
Further systematic computations would exceed reasonable time constraints.
The results obtained by our analysis give of course rise to further research.
It is still an open problem how graphs which have the maximal number of realizations can be classified,
and how to bound this number:
\begin{oproblem}
  Find an upper bound $b_n<4^{n-2}$ such that $\lamII(G)\leq b_n$ for all Laman graphs $G$ with $n$ vertices.
\end{oproblem}
From our data we observe that (for $n\leq12$) there is always a unique graph $G_{n,\max}$
on $n$ vertices that achieves the maximal Laman number among all graphs with $n$ vertices.
\begin{conjecture}
  For each $n\geq2$ there is a unique Laman graph $G_{n,\max}$ with $n$
  vertices and with the property $\lamII\left(G_{n,\max}\right)=\maxlamIIn$.
  Similarly for $\maxlamIIIn$.
\end{conjecture}
Also the relation of Henneberg steps to the increase of the number of realizations is subject of further research:
\begin{oproblem}
  Find lower and upper bounds for the factor $\lam(G')/\lam(G)$ where $G'$ is constructed from a Laman graph $G$ by a Henneberg step.
  By what we showed, the lower bound is smaller than or equal to $12/7$ in 2D and $1/4$ in 3D.
  The upper bound is bigger than or equal to $301/32$ in 2D and $68/3$ in 3D.
\end{oproblem}

In dimension 2 we expect every Henneberg step to increase the Laman number by at least a factor of two.
As mentioned above, this is still open for steps of type 2c.
\begin{conjecture}
  For a Laman graph $G$ with $n$ vertices we have $\lamII(G)\geq 2^{n-2}$.
\end{conjecture}
In dimension 3 this does definitely not hold any more, since the first line in Table~\ref{table:Laman_Increase_3d} gives a counterexample.
It would be interesting to know whether there is a lower bound on the Laman number in 3D.

Another direction of research is the study of real realizations, i.e., by
considering labelings~$\lambda$ whose values are in~$\R$ and embeddings into~$\R^d$.
In the 2D case, it is known that the ratio between the number of real and complex
realizations can be arbitrarily close to~$0$, by exhibiting a particular graph
(of $8$ vertices and Laman number~$90$), which provably cannot have as many real
realizations as complex ones~\cite{Jackson2018}, and by gluing this
graph arbitrarily often together.
\begin{oproblem}
  Let $R_2(G)$ be the maximal (finite) number of different real realizations in~$\R^2$
  of a Laman graph~$G=(V,E)$, that can be achieved for some real
  labeling~$\lambda\colon E\to\R$. Clearly, for a Laman graph~$G$ that is
  constructible by using only Henneberg steps of type~1, we have $R_2(G)=\lamII(G)$.
  But what can we say about the sequence $\bigl(\varphi_n\bigr){}_{n\geq2}$ of quotients
  \[
    \varphi_n := \frac{R_2\left(G_{n,\max}\right)}{\maxlamIIn},
  \]
  i.e., we are asking about the gap between real and complex realizations for
  graphs with maximal Laman number.  From \cite{EmirisMoroz} we know that
  $\varphi_n=1$ for $n\leq7$, but does $\varphi_n=1$ hold for all~$n$?
  Probably not. Do we have $\lim_{n\to\infty}\varphi_n=0$, or does this limit
  approach a nonzero constant?  Does the limit exist at all?
\end{oproblem}
Similar questions can be posed for the three-dimensional case, where
much less is known. A first step into this direction would be to answer
the following question:
\begin{oproblem}
  Find a Geiringer graph that cannot have as many real realizations as
  complex realizations.
\end{oproblem}

\minisec{Acknowledgments} 
G.\ Grasegger and C.\ Koutschan were partially supported by the Austrian Science Fund (FWF): W1214-N15, project DK9.
C.\ Koutschan was also supported by the Austrian Science Fund (FWF): F5011-N15.
E.~Tsigaridas is partially supported by ANR JCJC GALOP (ANR-17-CE40-0009).

\bibliographystyle{plain}
\bibliography{literature}
\newpage
\section*{Appendix --- Graph Encodings}\label{appendix:encoding}
In this section we present details on our graph encodings and collect the
encodings of the graphs explaining the results and observations in the main part.

We represent a graph by the integer that is obtained by flattening the upper right triangle
of its adjacency matrix and interpreting this binary sequence as an integer. Note
that the adjacency matrix will always have zeros on the main diagonal, and hence we
consider only entries above the main diagonal.
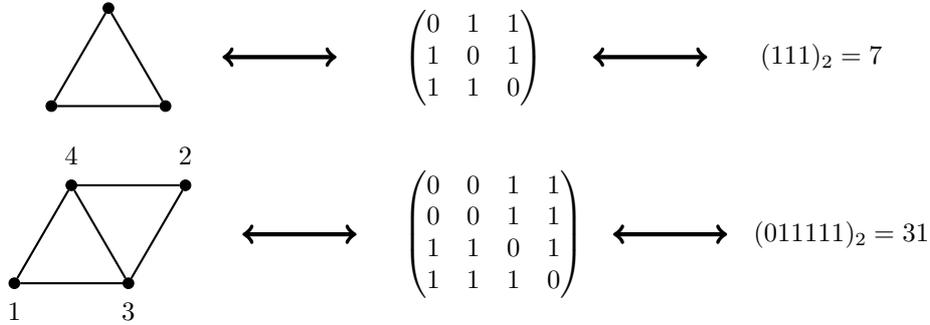
\begin{center}
  \begin{tikzpicture}[scale=1.5]
    \vertex (a) at (0,0) {};
		\vertex (b) at (1,0) {};
		\vertex (c) at (0.5,0.866025) {};
		\draw[edge] (a)edge(b) (a)edge(c) (b)edge(c);
		\draw[ultra thick,<->] (1.5,0.433) -- (2.5,0.433);
		\draw[ultra thick,<->] (4.75,0.433) -- (5.75,0.433);
		\node at (3.7,0.433) {
			$
			\begin{pmatrix}
			  0 & 1 & 1 \\
			  1 & 0 & 1 \\
			  1 & 1 & 0
			\end{pmatrix}
			$
		};
		\node at (6.75,0.433) {$(111)_2=7$};
  \end{tikzpicture}\\[2ex]
  \begin{tikzpicture}[scale=1.5]
    \vertex (a) at (-0.5,0) {};
    \node at (a.south) [below=2pt] {1};
		\vertex (b) at (0.5,0) {};
		\node at (b.south) [below=2pt] {3};
		\vertex (c) at (0,0.866025) {};
		\node at (c.north) [above=2pt] {4};
		\vertex (d) at (1,0.866025) {};
		\node at (d.north) [above=2pt] {2};
		\draw[edge] (a)edge(b) (a)edge(c) (b)edge(c) (b)edge(d) (c)edge(d);
		\draw[ultra thick,<->] (1.5,0.433) -- (2.5,0.433);
		\draw[ultra thick,<->] (4.75,0.433) -- (5.75,0.433);
		\node at (3.7,0.433) {
			$
			\begin{pmatrix}
			  0 & 0 & 1 & 1\\
			  0 & 0 & 1 & 1\\
			  1 & 1 & 0 & 1\\
			  1 & 1 & 1 & 0
			\end{pmatrix}
			$
		};
		\node at (6.75,0.433) {$(011111)_2=31$};
  \end{tikzpicture}
\end{center}
Note, that isomorphic graphs might be represented by different numbers in this way.
Hence, for our computations we used some normal form, which is not necessary to explain in detail here.
The conversion from a number to a graph does not depend on this normal form.

\begin{table}[ht]
  \begin{center}
    \begin{tabular}{lll}
      \toprule
      n  & Graph encoding     & Laman number\\\midrule
      6  & 7916               & 24 \\
      7  & 1269995            & 56 \\
      8  & 170989214          & 136 \\
      9  & 11177989553        & 344 \\
      10 & 4778440734593      & 880 \\
      11 & 18120782205838348  & 2288 \\
      12 & 252590061719913632 & 6180 \\
      \bottomrule
    \end{tabular}
  \end{center}
  \caption{Graph encodings for the graphs with maximal Laman number (see Figure~\ref{figure:max_graphs})}
  \label{table:enc:max_laman_numbers}
\end{table}
\begin{table}[ht]
  \begin{center}
    \begin{tabular}{lll}
      \toprule
      n  & Graph encoding                              & Laman number \\\midrule
      12 &  757486969329934592                         & 5952 \\
      13 & 3102079810848683155456                      & 15056 \\
      14 & 12393113433401056197689344                  & 39696 \\
      15 & 101535867160732294622504828928              & 105384 \\
      16 & 283980994531838217547205604229120           & 277864 \\
      17 & 65135173642079980743135145171586662400      & 731336 \\
      18 & 9061092056503516236392931137633162134437921 & 1953816 \\
      \bottomrule
    \end{tabular}
  \end{center}
  \caption{Graph encodings for the graphs in $T(n)$ (see Figure~\ref{figure:max_graphs_triangle})}
  \label{table:enc:max_graphs_triangle}
\end{table}

\begin{table}[ht]
  \begin{center}
    \begin{tabular}{lll}
      \toprule
      n  & Graph encoding & Laman number\\\midrule
      13 & 2731597771584836257824                      & 15536 \\
      14 & 3932631430916370534240769                   & 42780 \\
      15 & 94091005932357252120217796609               & 112752 \\
      16 & 892527555716690691964688718172672           & 312636 \\
      17 & 97035633928660816927022803757023440896      & 870414 \\
      18 & 1132478330239973528711451061872988363235584 & 2237312 \\
      \bottomrule
    \end{tabular}
  \end{center}
  \caption{Graph encodings for the graphs in $S(n)$ from Figure~\ref{figure:max_graphs_high}}
  \label{table:enc:max_graphs_high}
\end{table}

\begin{table}[ht]
  \begin{center}
    \begin{tabular}{lll}
      \toprule
      n  & Graph encoding                         & Laman number\\\midrule
      13 & 517844367551685511200                  & 15268   \\
      14 & 8465213527269428904345612              & 40088   \\
      17 & 34561064106536153162036856640676376576 & 1953816 \\
      \bottomrule
    \end{tabular}
  \end{center}
  \caption{Graph encodings for graphs which contain the triangle as subgraph and have high Laman number.}
  \label{table:enc:fan}
\end{table}

\begin{table}[ht]
  \begin{center}
    \begin{tabular}{lll}
      \toprule
      n  & Graph encoding                          & Laman number\\\midrule
      7  & 127575                                  & 48     \\
      8  & 7654183                                 & 112    \\
      9  & 11987422577                             & 288    \\
      10 & 26665598300033                          & 688    \\
      11 & 18226243755613920                       & 1760   \\
      12 & 57080320167818985484                    & 4864   \\
      13 & 1845359412452332949520                  & 12616  \\
      14 & 2116433716010931973523488               & 32984  \\
      15 & 366442648507105101448244891666          & 83792  \\
      16 & 1054776952932226148552313881544736      & 224976 \\
      17 & 260539761471154896904085679883542331426 & 570544 \\
      \bottomrule
    \end{tabular}
  \end{center}
  \caption{Graph encodings for the graphs from Figure~\ref{figure:sub4_max} and further graphs which contain the 4-vertex Laman graph as subgraph and have high Laman number.}
  \label{table:enc:31-fan}
\end{table}

\begin{table}[ht]
  \begin{center}
    \begin{tabular}{lll}
			\toprule
			n  & Graph encoding                          & Laman number\\\midrule
			6  & 3326                                    & 16 \\
			7  & 190686                                  & 32 \\
			8  & 210799326                               & 96 \\
			9  & 27047004894                             & 224 \\
			10 & 220302198846                            & 576 \\
			11 & 511412109882689                         & 1376 \\
			12 & 270814819769185025                      & 3648 \\
			13 & 2585030414085585133728                  & 9472 \\
			14 & 6356539347198988132306956               & 24752 \\
			15 & 1109200018557493535348018405392         & 62416 \\
			16 & 5598668013338146547621855406197248      & 168256 \\
			17 & 176789006904155934327358957938973624416 & 433920 \\
      \bottomrule
    \end{tabular}
  \end{center}
  \caption{Laman graphs which have the 5-vertex graph with encoding 254 as subgraph}
  \label{table:enc:254-fan}
\end{table}

\begin{table}[ht]
  \begin{center}
    \begin{tabular}{lllll}
			\toprule
			n  & Graph encoding       & Laman number & Graph encoding       & Laman number\\\midrule
			6  & 12511                & 16           & 10479                & 16 \\
			7  & 111335               & 32           & 103805               & 32 \\
			8  & 6419031              & 96           & 12339295             & 96 \\
			9  & 812960551            & 224          & 1024072271           & 224 \\
			10 & 209151514913         & 576          & 221350536519         & 576 \\
			11 & 110640260854593      & 1376         & 18441562579184833    & 1376 \\
			12 & 37616617704925531361 & 3648         & 21047011153048344071 & 3648 \\
      \bottomrule
    \end{tabular}
  \end{center}
  \caption{Laman graphs which have the 5-vertex graph with encoding 223 and 239 as subgraph, respectively}
  \label{table:enc:223-239-fan}
\end{table}

\begin{table}[ht]
  \begin{center}
    \begin{tabular}{lll}
			\toprule
			n  & Graph encoding                         & Laman number\\\midrule
			7  & 120478                                 & 48 \\
			8  & 6475132                                & 96 \\
			9  & 51946608057                            & 288 \\
			10 & 18284890201676                         & 672 \\
			11 & 5366995734673421                       & 1728 \\
			12 & 523614257391638273                     & 4128 \\
			13 & 2066305871268252766241                 & 10944 \\
			14 & 40197303758420411293510144             & 28416 \\
			15 & 61903368089062917457613881376          & 70656 \\
			16 & 11358585136343922383033065301099552    & 177408 \\
			17 & 33233417861308024077754506274593047824 & 486528 \\
      \bottomrule
    \end{tabular}
  \end{center}
  \caption{Laman graphs which have the three-prism with encoding 7916 as subgraph}
  \label{table:enc:7916-fan}
\end{table}

\begin{table}[ht]
  \begin{center}
    \begin{tabular}{lll}
      \toprule
      n  & Graph encoding & 3D-Laman number\\\midrule
      4  & 63             & 2 \\
      5  & 511            & 4 \\
      6  & 16350          & 16 \\
      7  & 515806         & 48 \\
      8  & 49724126       & 160 \\
      9  & 7345971057     & 640 \\
      10 & 3559487592083  & 2560 \\
      \bottomrule
    \end{tabular}
  \end{center}
  \caption{Graph encodings for the Geiringer graphs with maximal 3D-Laman number (see Figure~\ref{figure:max_graphs_3d})}
  \label{table:enc:max_graphs_3d}
\end{table}

\begin{table}[ht]
  \begin{center}
    \begin{tabular}{lll}
      \toprule
      n  & Graph encoding & 3D-Laman number\\\midrule
      5  & 511            & 4 \\
      6  & 7679           & 8 \\
      7  & 257911         & 32 \\
      8  & 16559991       & 96 \\
      9  & 4076665507     & 448 \\
      10 & 4894450217603  & 1664 \\
      \bottomrule
    \end{tabular}
  \end{center}
  \caption{Graph encodings for the Geiringer graphs which contain the tetahedron}
  \label{table:enc:fan_3d}
\end{table}

\begin{table}[ht]
  \begin{center}
    \begin{tabular}{lll}
      \toprule
      n  & Graph encoding & 3D-Laman number\\\midrule
      6  & 7679           & 8 \\
      7  & 237055         & 16 \\
      8  & 14937975       & 64 \\
      9  & 38164887119    & 256 \\
      10 & 3168405805643  & 896 \\
      \bottomrule
    \end{tabular}
  \end{center}
  \caption{Graph encodings for the Geiringer graphs which contain the double tetahedron}
  \label{table:enc:genfan_3d}
\end{table}

\end{document}